\documentclass[letterpaper, 10 pt, conference]{ieeeconf}  

\IEEEoverridecommandlockouts                              

\usepackage{balance}
\usepackage{cite}
\usepackage{amsmath,amsfonts}

\usepackage{todonotes}

\newtheorem{example}{Example}

\title{\LARGE \bf
Direct Transcription for Dynamic Optimization: A Tutorial with a Case Study on Dual-Patient Ventilation During the COVID-19 Pandemic
}

\author{
Eric C. Kerrigan$^{1,2}$, Yuanbo Nie$^{1}$, Omar Faqir$^{2}$, Caroline H. Kennedy$^{3}$, 
\\
Steven A. Niederer$^{4}$, Jose A. Solis-Lemus$^{4}$, Peter Vincent$^{1}$, Steven E. Williams$^{4,5}$
\thanks{This work was supported by the EPSRC Centre for Doctoral Training in High Performance Embedded and Distributed Systems  (HiPEDS, Grant Reference EP/L016796/1) and the Wellcome/EPSRC Centre for Medical Engineering (WT 203148/Z/16/Z). }
\thanks{$^{1}$E. C. Kerrigan, Y. Nie  and P. Vincent are with the Department of Aeronautics, Imperial College London, UK,  e.kerrigan@imperial.ac.uk, yuanbo.nie15@imperial.ac.uk, p.vincent@imperial.ac.uk}
\thanks{$^{2}$E. C. Kerrigan and O. Faqir are with the Department of Electrical \& Electronic Engineering, Imperial College London, UK,  omar.faqir12@imperial.ac.uk}
\thanks{$^{3}$C. H. Kennedy is with the Evelina Children’s Hospital, Guy’s and St Thomas’ NHS Foundation Trust, London, UK, Caroline.Kennedy@gstt.nhs.uk}
\thanks{$^{4}$S. A. Niederer, J. A. Solis-Lemus and S. E. Williams are with the School of Biomedical Engineering and Imaging Sciences, King’s College London, UK,  steven.niederer@kcl.ac.uk, jose.solislemus@kcl.ac.uk, steven.e.williams@kcl.ac.uk}
\thanks{$^{5}$S. E. Williams is also with the Department of Cardiology, Guy’s and St Thomas’ NHS Foundation Trust, London,  UK}
}

\begin{document}

\sloppy
\maketitle
\thispagestyle{empty}
\pagestyle{empty}




\begin{abstract}
   A variety of optimal control, estimation, system identification and design problems can be formulated as functional optimization problems with differential equality and inequality constraints. Since these problems are infinite-dimensional and often do not have a known analytical solution, one  has to resort to numerical methods to compute an approximate solution. This paper uses a unifying notation to outline some of the techniques used in the  transcription step of simultaneous direct  methods (which discretize-then-optimize) for solving continuous-time dynamic optimization problems. We focus on collocation, integrated residual and Runge-Kutta schemes. These transcription methods are then applied to a simulation case study to answer a question  that arose during the COVID-19 pandemic, namely: If there are not enough ventilators, is it possible to ventilate more than one patient on a single ventilator? The results suggest that it is possible, in principle, to estimate individual patient parameters sufficiently accurately, using a relatively small number of flow rate measurements, without needing to disconnect a patient from the system or needing more than one flow rate sensor.  We also show that it is possible to ensure that two different patients can indeed receive their desired tidal volume, by modifying the resistance experienced by the air flow to each patient and controlling the ventilator pressure.
\end{abstract}

\section{Introduction}

\subsection{Problem Formulation}

Many optimal control, estimation, 
system identification and system design problems can be formulated as a finite-horizon  \emph{dynamic optimization problem} (DOP). That is, one seeks to optimize and constrain the evolution of a dynamical system on a time interval $\mathcal{T}:=[t_0,t_f]\subset \mathbb{R}$, where $t_0$ and~$t_f$ denote the initial and final time, respectively.

We will consider continuous-time DOPs that can be written in the popular Bolza form~\cite{betts2010practical,RMD2nd}, i.e.\ we seek to find solutions to the optimization problem
\begin{subequations}
\label{eqn:cont_DOP}
\begin{equation}
    \min_{\substack{x(\cdot),u(\cdot)\ \\ \theta,t_0,t_f}}
    V_M(x(t_0),x(t_f),\theta,t_0,t_f) +  
    \int_{t_0}^{t_f}\ell(x(t),u(t),\theta,t) dt
\end{equation}
subject to the state trajectory $x: \mathbb{R} \rightarrow \mathbb{R}^{n_x}$ being continuous and the following constraints being satisfied:
\begin{align}
    f(\dot{x}(t),x(t),u(t),\theta,t) = 0,\ &\forall t \in\mathcal{T}\ \text{a.e.}, \label{eqn:DEs}\\
    g(\dot{x}(t),x(t),\dot{u}(t),u(t),\theta,t) \leq 0,\  &\forall t \in\mathcal{T} \ \text{a.e.}, \label{eqn:ineqs}\\
    c(\dot{x}(t),x(t),\dot{u}(t),u(t),\theta,t) = 0,\  &\forall t \in\mathbb{T}^c,  \label{eqn:wp_const}\\
    \psi_E(x(t_0),x(t_f),\theta,t_0,t_f) = 0,\ & \label{eqn:terminal_eq}\\
    \psi_I(x(t_0),x(t_f),\theta,t_0,t_f)\leq 0,\ \label{eqn:terminal_ineq}&
\end{align}
\end{subequations} 
where  `a.e.' stands for `almost everywhere' in the Lebesgue sense. In other words, the state trajectory $x$ and trajectory of free variables $u: \mathbb{R} \rightarrow \mathbb{R}^{n_u}$ is allowed to be non-differentiable (but is piecewise differentiable) and~\eqref{eqn:DEs}--\eqref{eqn:ineqs} is allowed to be violated on a set of measure zero.  

In control and system design problems, the vector $u(t)$ includes so-called control inputs or manipulated variables, which are time-varying physical or virtual variables that can be adjusted by a human or automatic control system. Examples include the amount of power or fuel used, rate at which money is spent, percentage of a population quarantined, actuator position or force applied at time $t$. The function $u$ can also include time-varying parameters used to define feedforward and feedback policies in 
robust and stochastic optimal control problems~\cite{HandbookMPC}. In estimation and system identification problems $u$ is often used to model external (sometimes called uncontrolled) inputs, such as unknown disturbances, measurement noise or unknown time-varying parameters. 

The vector $\theta \in \mathbb{R}^{n_s}$ includes all constant parameters to be determined. For example, in control problems $\theta$ could include the amount of energy at the start, capital to invest or parameters of feedforward and feedback policies to be determined in multiple-scenario problems. In estimation and system identification problems~$\theta$  could include the parameters in linear and nonlinear black- or grey-box models. In system design problems~$\theta$ could include parameters to be determined during design time, such as the size of a battery or the geometry and mass of an object.

The evolution of the system is assumed to be described by ordinary differential equations (ODEs) or differential-algebraic equations (DAEs) that can be written in the form~\eqref{eqn:DEs}, where $f:\mathbb{R}^{n_x} \times \mathbb{R}^{n_x} \times \mathbb{R}^{n_u} \times \mathbb{R}^{n_s} \times \mathbb{R} \to \mathbb{R}^{n_f}$. To simplify notation,  algebraic variables for DAEs are assumed to be included in the free variables $u$, rather than introducing another variable to represent them.  Recall also that in many applications the DAEs or ODEs arise from the semi-discretization of partial differential equations (PDEs); how best this should be done is outside the scope of this paper, but some of the discussion here will still be applicable.

The cost function includes a so-called Mayer cost term $V_M: \mathbb{R}^{n_x} \times \mathbb{R}^{n_x} \times \mathbb{R}^{n_s} \times \mathbb{R} \times \mathbb{R} \to \mathbb{R}$ and the integral of the running cost $\ell:\mathbb{R}^{n_x} \times \mathbb{R}^{n_u} \times \mathbb{R}^{n_s} \times \mathbb{R} \to \mathbb{R}$. These terms can be used to model a variety of popular functions to minimise. Examples of suitable cost functions in control and system design include time spent to complete a task, energy used, size of infected population or money spent. Examples of suitable cost functions in estimation and system identification include weighted least-squares terms for the noise and disturbance that explains the mismatch between the measurements and the model.


In control and  design problems the system is often subject to additional constraints, e.g.\ upper and lower bounds on the actuators or constraints arising due to performance, legal or safety specifications. In estimation and system identification problems it is also common practice to assume that the disturbances, noise or parameters satisfy certain constraints. We assume that all  inequality constraints are  captured in~\eqref{eqn:ineqs}, where $g:\mathbb{R}^{n_x} \times \mathbb{R}^{n_x} \times \mathbb{R}^{n_u} \times \mathbb{R}^{n_u} \times \mathbb{R}^{n_s} \times \mathbb{R} \to \mathbb{R}^{n_g}$.

In some control applications there is  a finite set of constraints given as  equality constraints, e.g.\ where a robotic end effector has to pass through a certain sequence of points in space. In estimation and system identification problems, a  finite set of noisy measurements is usually given. This is captured with~\eqref{eqn:wp_const}, where  $\mathbb{T}^c$  is a finite subset of $\mathcal{T}$, and $c:\mathbb{R}^{n_x} \times \mathbb{R}^{n_x} \times \mathbb{R}^{n_u} \times \mathbb{R}^{n_u} \times \mathbb{R}^{n_s} \times \mathbb{R} \to \mathbb{R}^{n_c}$.

 We also   include separate boundary equality and inequality constraints, as well as constraints on the initial and final time, in~\eqref{eqn:terminal_eq}--\eqref{eqn:ineqs}, where
$\psi_E:\mathbb{R}^{n_x} \times \mathbb{R}^{n_x} \times \mathbb{R}^{n_s} \times \mathbb{R} \times \mathbb{R} \to \mathbb{R}^{n_E}$ and $\psi_I:\mathbb{R}^{n_x} \times \mathbb{R}^{n_x} \times \mathbb{R}^{n_s} \times \mathbb{R} \times \mathbb{R} \to \mathbb{R}^{n_I}$. 

To simplify notation, we focus on problems where~$x$ is constrained to be continuous over the whole interval $\mathcal{T}$. See~\cite{betts2010practical} on how to handle so-called multi-phase or hybrid problems, where a discontinuous trajectory $x$ is allowed, e.g.\ in multi-stage rockets or walking robots. 
Note  that it is possible  to add a constraint that $u$, $\dot{x}$ or $\dot{u}$ be continuous by adding constraints to the discretized problem in the same manner as which constraints are added to ensure continuity of the state (discussed in Section~\ref{sec:parameterization_x}).

\subsection{Scope, Aims and Case Study}

It can be very challenging to solve DOPs in the above form. The optimization problem is infinite-dimensional, because we are seeking to optimize over functions $x$ and $u$ that live in an infinite-dimensional space and there is an uncountable set of  constraints. Hence, analytical solutions often do not exist for practical problems. In these cases, sometimes the only  way forward is to use numerical methods to compute approximate solutions. This is the topic of the paper. Section~\ref{sec:transcription} will focus on describing some of the most popular classes of methods that allow a designer to directly discretize the above problem in order to compute an approximate solution using numerical optimization methods. 

The main aim of this paper is to introduce a unifying, abstract framework and notation in which a selection of direct transcription methods could be introduced to a non-expert. 
The paper is tutorial in nature and is not intended to be a comprehensive survey or review.
Some implementation details have unfortunately had to be  omitted in order to prevent the paper from turning into a book.
However,  we hope that the presentation allows both the non-expert and  expert alike to come to a clearer understanding as to what the key concepts of a method are, as well as what the similarities and differences between some methods are. 

A reader of a tutorial paper expects either one challenging or many small example problems. At the time of writing, we found ourselves in the first few months of the COVID-19 pandemic. We were faced with a question that many scientists and engineers around the world were asking: What can a hospital do if there are not enough ventilators for patients? Because there were a number of unexpected demands on our time and we wanted to contribute to finding an answer to this question, we decided to focus our efforts in this tutorial paper on  presenting a mathematical description of the application we were working on at the time.
Sections~\ref{sec:ventilators} and~\ref{sec:results} therefore extends some initial simulation results of ours from a multi-disciplinary study aimed at  developing a system for ventilating more than one patient on a single ventilator~\cite{doi:10.1098/rsos.200585}. 
The study aims to understand how best to monitor and manipulate the air flow to individual patients when they are attached to the same ventilator.

This problem has  the following features, which motivate the application of state-of-the-art methods in   system identification, estimation and control:
\begin{itemize}
    \item The dynamics are nonlinear and it is  an open question whether linear models are sufficiently accurate.
    \item Model parameters are unknown and there are neglected dynamics.
    \item Humans, namely both patients and clinical staff, are part of the closed-loop system.
    \item A wide range of scenarios and solutions may be considered, ranging from ones that require minimum changes to current ventilator setups, to future ventilator designs that could have a splitting option enabled by the manufacturer.
\end{itemize}

We would like to stress that the focus of this paper is as a tutorial paper on transcription methods for solving dynamic optimization problems. The case study is not meant to represent the state-of-the-art, a review or survey on patient ventilation. Experts on ventilator design are likely to be disappointed. The case study was chosen  as a topical and challenging problem, which we were learning about ourselves while writing this paper, yet was amenable to treatment in a tutorial context. Our hope is that the case study will also  initiate some  new scientific, engineering, medical and ethical questions, to which we do not  have  answers.

\section{Direct Transcription Methods}
\label{sec:transcription}

Direct transcription is the procedure whereby the continuous-time dynamic optimization problem~\eqref{eqn:cont_DOP}, which is an infinite-dimensional optimization problem, can be  `approximated' by a 
finite-dimensional 
optimization problem of the form:
\begin{subequations}
\label{eqn:dt_DOP}
\begin{align}
    \min_{q,s,\pi}\
    V_M^d(s_0,s_N,\pi)+
    \sum_{i=0}^{N-1}L_i(s_i,q_i,\pi)
    \label{eqn:dt_DOP_cost}
\end{align}
subject to
\begin{align}
    \phi_i(s_i,s_{i+1},q_i,\pi) &= 0,\ \forall i \in \mathbb{I}_N, \label{eqn:dt_DOP_eq}\\
    \gamma_i(s_i,q_i,\pi) &\leq 0,\ \forall i \in \mathbb{I}_N,\ &\label{eqn:dt_DOP_ineq}\\
     \psi_E^d(s_0,s_N,\pi) &= 0,\\
    \psi_I^d(s_0,s_N,\pi)&\leq 0, \label{eqn:dt_DOP_ineq_terminal}
\end{align}
\end{subequations}
where\footnote{Given  column vectors $a$ and $b$, the notation $(a,b):=[a^\top\ b^\top]^\top$.}   $\pi:=(\theta,t_0,t_f)$, $s_i\in\mathbb{R}^{n_i}$, $q_i\in\mathbb{R}^{m_i}$,   for all  $i \in \mathbb{I}_N:=\{0,\ldots,N-1\}$,  $s:=(s_0,\ldots,s_N)$ and $q:=(q_0,\ldots,q_{N-1})$.
The functions $V_M^d$, $L_i$, $\phi_i$, $\gamma_i$, $\psi_E^d$ and $\psi_I^d$ are discretized forms of~$V_M$, $\ell$, $f$, $g$, $\psi_E$ and $\psi_I$. 
The above problem is highly structured and can be efficiently solved using state-of-the-art nonlinear programming (NLP) solvers~\cite{FERREAU201713194}, which will be  discussed in Section~\ref{sec:NLP}.

Note that the lengths of the vectors $q_i$ and $s_i$ are  functions of the stage number $i$. 
A convenient and popular choice for  $q_i$ and $s_i$ are for these vectors to be composed of sampled versions of the  trajectory of free variables~$u$ and the state trajectory  $x$, respectively. However, this is not the only possibility --- the definitions of $q$ and $s$ depend on the  parameterization of the trajectories in the  transcription procedure, as will be discussed below. 

Note also that $\phi_i$ an $\gamma_i$ are not functions of $f$ and $g$, respectively. As will be seen below,  they could be functions of both $f$ and $g$, or neither.

This paper is mostly concerned with providing an introduction to a class of so-called simultaneous direct transcription methods, which  translates the continuous-time problem into a single, stand-alone NLP.
For completeness, we also discuss direct shooting methods in Section~\ref{sec:shooting}. These methods differ from the simultaneous schemes considered in Sections~\ref{sec:collocation} to~\ref{sec:RK} in that they solve the differential equations by interfacing the NLP  solver to separate, stand-alone numerical differential equation solvers. In contrast, the class of methods we focus on here  solve the differential equations without the use of separate differential equation solvers.

\subsection{Variable-time Problems}

In variable-time problems, such as minimum-time problems, the final time $t_f$ and/or initial time $t_0$ are decision variables, subject to given inequality constraints on them. There are a number of different ways to handle this when transcribing the problem into a finite-dimensional optimization problem. One popular way is to employ a transformation on the time so that the problem becomes one with fixed start and end times. 
For this purpose, let  each time instance $t\in\mathcal{T}$ be associated with a fixed, non-dimensional time instance $\tilde{t} \in [0,1]$ such that
\[
t = t_0+\tilde{t}(t_f-t_0) \Leftrightarrow \tilde{t} = (t-t_0)/(t_f-t_0).
\]
If the above substitution of variables is made in~\eqref{eqn:cont_DOP}, then the new problem becomes one with $\tilde{t}$ as the new time variable, with fixed start time $\tilde{t}_0=0$ and fixed end time~$\tilde{t}_f=1$. The variables $t_0$ and $t_f$ are then included as part of the parameters~$\pi$ in the discretized DOP~\eqref{eqn:dt_DOP}.

In order to simplify the presentation in this paper, but pointing out to the reader that variable-time problems can be handled in the schemes presented below,  the explicit dependence of a time instance~$t_i$ on  $\tilde{t}_i$, $t_0$ and $t_f$ will be omitted. It should  be understood throughout that, if the problem is variable-time, then the time instances in any given finite subset of $\mathcal{T}$, and therefore the differences between time instances, are actually functions of $t_0$ and~$t_f$; in contrast, each~$\tilde{t}_i$ and the differences between them are not functions of $t_0$ and~$t_f$.

\subsection{Approximate Parameterization of $u$}
\label{sec:parameterization_u}
\begin{figure}[t]
\begin{center}
\includegraphics[width=0.7\columnwidth]{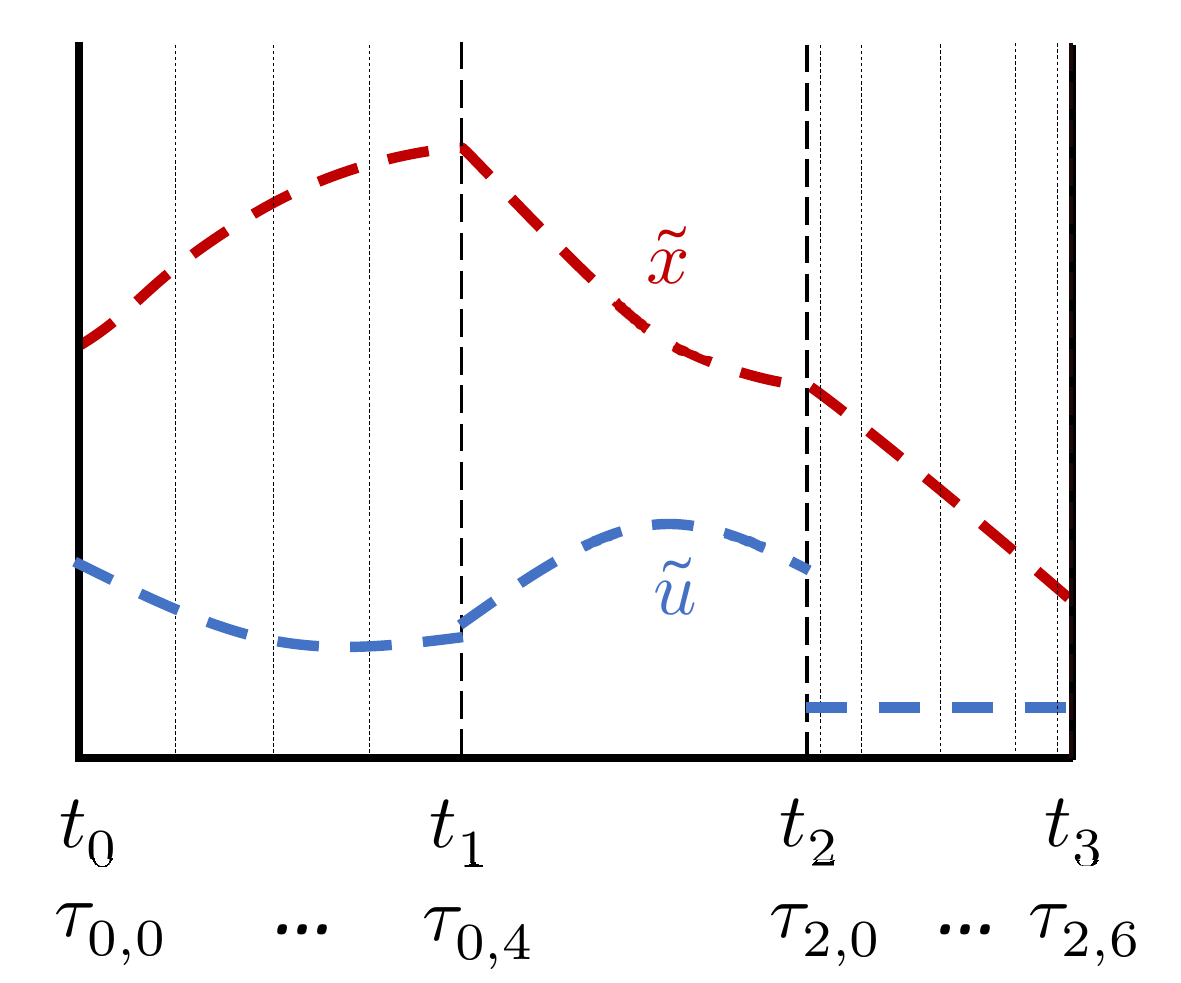}
\caption{Temporal discretization and trajectory parameterization. Note that the number and distribution of interpolation points need not be the uniform or the same in each interval. Note that the state trajectory $\tilde{x}$ is continuous, but that the trajectory $\tilde{u}$ can be discontinuous.} 
\label{fig:ParameterizationAndMesh}
\end{center}
\end{figure}

A sensible first step in a transcription process is to parameterize a finite-dimensional approximation of the  trajectory of free variables $u$ via the components of~$q$.  The discretization and parameterization for $u$ and $x$ is illustrated in Figure~\ref{fig:ParameterizationAndMesh}.
The development here is a generalization of what is  referred to as \emph{move blocking} in the predictive control literature~\cite{iet:/content/journals/10.1049/iet-cta.2019.0168}, where the input trajectory is typically constrained to be piecewise constant on intervals of different lengths.

We start by defining the set of mesh nodes
\[
    \mathbb{M}:= \{t_0,t_1,t_2,\ldots,t_N\},
\]
where
\[
    t_0< t_1 < \cdots < t_N = t_f
\]
and the length of each interval $\mathcal{T}_i:=[t_i,t_{i+1}]$ is
\[
h_i := t_{i+1}-t_i,\ \forall i \in\mathbb{I}_N.
\]

The trajectory of free variables $u$ is  approximated by
\begin{equation}
u(t)\approx \tilde{u}(t;q) := \upsilon_i(t;q_i),\ \forall t\in[t_i,t_{i+1}), i\in\mathbb{I}_N,
\end{equation}
where each function $\upsilon_i$ is  continuous and 
differentiable; suitable choices are discussed below. Recall that $u$ can be discontinuous; w.l.o.g.\ we defined $\tilde{u}$ to take on the left-hand limit at $t_i$ in the above, as is  convention. However, note that the left limit $\upsilon_i(t_{i+1};q_i)$ at $t_{i+1}$, and not the right limit $\upsilon_{i+1}(t_{i+1};q_{i+1})$, should be used when  transcribing the problem in the interval $\mathcal{T}_i$; this is particularly important when using implicit Runge-Kutta, Radau or Lobatto schemes, such as the trapezoidal and Hermite-Simpson methods, all of which evaluate functions at the boundary of the interval.

Suitable choices for the functions $\upsilon_i$ are application-dependent. Popular choices are for each $t \mapsto\upsilon_i(t;q_i)$ to be a constant or polynomial, so that the parameterized trajectory $t\mapsto \tilde{u}(t;q)$ is   piecewise constant or (discontinuous) piecewise  polynomial. However, other parameterizations could also be used, e.g.\ piecewise trigonometric polynomial,  piecewise algebraic, piecewise transcendental or any combination of the above, if there is reason to justify that this might result in a more efficient or numerically reliable scheme than piecewise polynomials. In other words, $q_i$ could consist of the coefficients associated with a suitable set of basis functions or the weights and biases of a neural network.

Note that in some control systems it might be necessary for $u$ to be implemented with a zero-order or other hold. However, early on in the mesh refinement process (see Section~\ref{subsec:Error_And_MeshRefinement}), when the intervals defined by the mesh are relatively large compared to the time between control updates, it might be computationally more efficient to use a  (non-constant) polynomial parameterization for $\upsilon_i$,  even if $u$ is implemented in a piecewise constant manner.

A reasonable constraint to impose on the choice of parameterization is to require that~$q_i$ consists of samples of $\upsilon_i(\cdot;q_i)$, i.e.\ the trajectory $t\mapsto\upsilon_i(t;q_i)$ should interpolate through the components of
\[
q_i=:(q_{i,0},q_{i,1},\ldots,q_{i,N_i^q})\in\mathbb{R}^{m_i}\] 
at a given set of nodes
\[
\mathbb{T}_i^q=:\{\tau_{i,0}, \ldots, \tau_{i,{N_i^q}}\}\subset\mathcal{T}_i,
\]
where each vector $q_{i,j}\in\mathbb{R}^{n_u}$, 
so that $m_i = (N_i^q+1)n_u$ and
\[
    \upsilon_i(\tau_{i,j};q_i) = q_{i,j},\ \forall i\in\mathbb{I}_N,j \in\mathbb{J}_i^q:=\{0,\ldots, N_i^q\}.
\]
This is possible with Lagrange polynomials (see Appendix~\ref{app:Lagrange}).
One could therefore define $\upsilon_i(\cdot;q_i)$   as the Lagrange interpolating polynomial
\[
\upsilon_i(t;q_i):=
\sum_{j\in\mathbb{J}_i}q_{i,j} \mathcal{L}_j(t) ,\
\forall t\in\mathcal{T}_i,
\]
where the Lagrange polynomial
\[
\mathcal{L}_j(t):=\prod_{k \in\mathbb{J}_i^q,k\neq j} \frac{t-\tau_{i,k}}{\tau_{i,j}-\tau_{i,k}}, \ \forall j \in\mathbb{J}_i^q,t\in\mathcal{T}_i,
\]
has degree at most $N_i^q$. For $N_i^q=0,1,2,3,\ldots$ the function $t\mapsto \upsilon_i(t;q_i)$ is therefore constant, affine, quadratic, cubic, etc. 
As mentioned above, other choices for basis functions are available and it is not necessary for components of~$q_i$ to interpolate $\upsilon_i(\cdot,q_i)$. However,  warm starting of the NLP solver is significantly faster and simpler to implement with interpolating polynomials, compared to using non-interpolating parameterizations, such as monomial basis functions. Furthermore, Lagrange polynomials have a number of advantages from a numerical and implementation point of view, compared to using monomials or other functions as basis functions~\cite{doi:10.1137/S0036144502417715}, hence Lagrange polynomials are often the basis functions of choice. 

\subsection{Approximate Parameterization of $x$}
\label{sec:parameterization_x}
We proceed in a similar manner as in the previous section, but with some slight differences. In some cases it might be a good idea to use a different mesh for the states. However, for this tutorial we will simplify notation considerably by adopting the same set of mesh nodes $\mathbb{M}$ as above. 

We seek to approximate the state trajectory $x$  as
\begin{subequations}
\begin{equation}
    x(t)\approx \tilde{x}(t;s) := \chi_i(t;s_i),\ \forall t\in [t_i,t_{i+1}), i\in\mathbb{I}_N.\\
\end{equation}
where the functions $\chi_i$ are continuous and  differentiable.
We  also  define 
\begin{equation}
    \tilde{x}(t_f;s) :=s_N.
\end{equation}
\end{subequations}


In order to ensure that the trajectory $t\mapsto\tilde{x}(t;s)$ is continuous at the boundaries of each interval~$[t_i,t_{i+1})$ and at the final time, the following constraints have to be included in the equality constraints~\eqref{eqn:dt_DOP_eq}:
\begin{subequations} 
\label{eqn:continuity}
\begin{align}
\chi_{i-1}(t_i;s_{i-1}) &= \chi_i(t_i;s_i),\ \forall i\in\mathbb{I}_N\backslash \{0\}, \label{eqn:cont_i} \\
 \chi_{N-1}(t_f;s_{N-1}) &= s_N  \label{eqn:cont_N}.
\end{align}
\end{subequations}
As can be seen, the addition of the above  constraints is the reason for the coupling between  stages $i$ and~$i+1$ in~\eqref{eqn:dt_DOP_eq}.


As in the previous section, one could parameterize each $\chi_i$ such that the trajectory $t\mapsto\chi_i(t;s_i)$ is polynomial/trigonometric/algebraic/transcendental/other. The trajectory $t\mapsto\tilde{x}(t;s)$ will then be continuous and piecewise polynomial/trigonometric/algebraic/transcendental/other, as illustrated in Figure~\ref{fig:ParameterizationAndMesh}. 

It is also a good idea to choose the parameterization such that the trajectory $t\mapsto\tilde{x}(t;s)$ interpolates through some of the components of $s$ on a given set of nodes $\mathbb{T}_i^s$, e.g.\ using Lagrange interpolating polynomials to define some of the components of~$\chi_i$. However, as discussed in Section~\ref{sec:RK}, sometimes it is not necessary, possible or desirable that $\upsilon_i(\cdot;q_i)$ and $\chi_i(\cdot;s_i)$ interpolate through all components of~$q_i$ or $s_i$.

\subsection{Error Analysis and Mesh Refinement} \label{subsec:Error_And_MeshRefinement}

Given a solution to~\eqref{eqn:dt_DOP}, the level of accuracy and constraint satisfaction of the solution  needs to be checked. The size of the violation of the differential equations over each interval in-between grid points can be computed with the integral
\begin{equation*}
\zeta_i:= \frac{1}{h_i} \int_{\mathcal{T}_i}	\| f(\dot{\chi_i}(t;s),\chi_i(t;s),\upsilon_i(t;q),\theta,t)\|_p\: dt,
\end{equation*}
for all $i\in\mathbb{I}_N$. Alternatively, one could compute
\begin{equation*}
\zeta_{i,j}:=\frac{1}{h_i}\int_{\mathcal{T}_i}  |\ f_j(\dot{\chi_i}(t;s),\chi_i(t;s),\upsilon_i(t;q),\theta,t)|\: dt,
\end{equation*}
for each   component $j=1,\hdots,n_f$ of the  \emph{residual}~$f$. Note that the above integrals can be computed exactly or using very efficient quadrature schemes.

In either form, $\zeta$ is typically referred to as the \emph{absolute local error}. Another criterion commonly used is the \emph{relative local error} based on a selection of normalization weights (see~\cite{betts2010practical} for  details). Additionally, \emph{inequality constraint violation errors} can be computed to measure the extent of possible inequality constraint violations of the trajectories in-between the  points where the inequality constraints are enforced. 

When formulating the DOP, the  practitioner may indicate desired error levels, based on the above-mentioned criteria, in the form of error tolerances. Once the magnitudes and distributions of errors are identified, corresponding modifications to the grid can be made. The problem is iteratively solved until the solution fulfils all error tolerances. This process is called mesh refinement, generally requiring mesh nodes to be added to $\mathbb{M}$. One can also choose to increase the number of parameters/basis functions for the parameterizations for some or all of the $\chi_i$ and/or $\upsilon_i$ (see \cite{betts2010practical,liu2017adaptive} for example). The NLP problem based on the new mesh and parameterization can be warm-started using the solution from the previous mesh, leading to significantly faster convergence to a solution, thus reducing the overall computation time. The mesh is considered sufficiently fine if the above integrals are below specified tolerances and the magnitude of the difference between the cost with one mesh and a finer mesh is below a tolerance.

Another way to reduce the computation time under the mesh refinement framework is to use an external constraint handling scheme \cite{ECH_CSL:2020,polak:2009}. These schemes systematically include or exclude inactive inequality constraints that do not contribute to the optimal solution, but burden the numerical computations. 

It is important to stress that certain methods of error analysis typically employed when solving differential equations do not work for dynamic optimization problems, in general. For example, if the solutions to~\eqref{eqn:cont_DOP} and/or~\eqref{eqn:dt_DOP} are non-unique, then one cannot compare the solutions obtained using two different meshes to assess the accuracy of the solution. If the solution of the differential equation is unstable, as is often the case, then choosing to  keep the trajectory of external/control input variables fixed and solving for the trajectory of state and algebraic variables on a finer mesh, can fail to provide numerically reliable error estimates. Hence, the above integrals are commonly used for error analysis, instead of the usual error analysis methods employed when solving differential equations. 

\subsection{Collocation Schemes}
\label{sec:collocation}

Collocation is arguably the easiest direct transcription method to implement. Explicit/forward and implicit/backward Euler schemes are probably the best-known collocation methods. Other well-known methods include the mid-point rule, trapezoidal, Hermite-Simpson and pseudo-spectral LG, LGR and LGL methods.

In collocation methods, the trajectories $\tilde{u}(\cdot;q)$ and $\tilde{x}(\cdot;s)$ are required to satisfy the equality constraints~\eqref{eqn:DEs} exactly at only a finite subset of each interval~$\mathcal{T}_i$, namely at a set of so-called  \emph{collocation points}
\[
    \mathbb{T}_i^f:=\{
    t_i+h_i(1+c)/2 
    \mid c\in \mathbb{K}_i\}, \forall i \in \mathbb{I}_N,
\]
where each $\mathbb{K}_i$ is a given finite subset of the interval~$[-1,1]$ so that the coefficient $(1+c)/2\in[0,1]$.
A Radau scheme includes only one of the boundaries of the interval $[-1,1]$ in~$\mathbb{K}_i$. A Gauss scheme does not include any of the boundary points of $[-1,1]$ and a Lobatto scheme includes both $-1$ and~$1$. Particular types of polynomial collocation methods are discussed in more detail at the end of this section. 

The equality constraints~\eqref{eqn:dt_DOP_eq} are given by combining the  finite set of constraints
\begin{subequations}
\label{eqn:collocation}
\begin{multline}
 f(\dot{\chi}_i(t;s_i),\chi_i(t;s_i),\upsilon_i(t;q_i), \theta,t) = 0,\\ \forall t \in\mathbb{T}_i^f,i\in\mathbb{I}_N
\label{eqn:collocation_f}
\end{multline}
with the continuity constraints~\eqref{eqn:continuity}
and 
\begin{multline}
    c(\dot{\chi}_i(t;s_i),\chi_i(t;s_i),\dot{\upsilon}_i(t;q_i),\upsilon_i(t;q_i),\theta,t) 
    =0,\\ \forall t\in\mathbb{T}^c\cap\mathcal{T}_i,i\in\mathbb{I}_N.
    \label{eqn:wp_collocation}
\end{multline}
In some applications it might be convenient to use interpolating parameterizations and set the interpolation points to be the same as the set of points at which the equality constraints are enforced, in which case $\mathbb{T}_i^q=\mathbb{T}_i^s=\mathbb{T}_i^f\cup(\mathbb{T}^c\cap\mathcal{T}_i)\cup\{t_i\}$; however, this is not necessary. Care has to be taken, though, to ensure that the system of equations is not over-determined and that there are sufficient degrees of freedom for a solution to exist, which will be discussed below. 

To enforce the inequality constraints~\eqref{eqn:ineqs}, let~\eqref{eqn:dt_DOP_ineq} be given by the finite set of constraints
\begin{multline}
    g(\dot{\chi}_i(t;s_i),\chi_i(t;s_i),\dot{\upsilon}_i(t;q_i),\upsilon_i(t;q_i),\theta,t) 
    \leq  h_i\varepsilon_i^g,\\ \forall t\in\mathbb{T}_i^g,i\in\mathbb{I}_N
    \label{eqn:ineq_collocation}
\end{multline}
where $\mathbb{T}_i^g$ is \emph{any} finite subset of~$\mathcal{T}_i$.  A convenient choice   is to enforce the inequality constraints only at the collocation points, i.e.\ $\mathbb{T}_i^g=\mathbb{T}_i^f$, but this is neither necessary nor sufficient to guarantee that the inequality constraints are satisfied in-between collocation points.
The vector $\varepsilon_i^g\leq 0$ therefore serves the purpose of tightening some or all of the constraints; a well-chosen $\varepsilon_i^g$ and $\mathbb{T}_i^g$ ensures that the inequality constraints~\eqref{eqn:ineqs} are satisfied at all time instances not in~$\mathbb{T}_i^g$~\cite{FP2019}.

The remaining constraints in~\eqref{eqn:dt_DOP} are given by
\begin{align}
    \psi_E^d(s_0,s_N,\pi):=\psi_E(\chi_0(t_0;s_0),s_N,\pi) &= 0,\label{eqn:collocation_boundary}\\
    \psi_I^d(s_0,s_N,\pi):=\psi_I(\chi_0(t_0;s_0),s_N,\pi) &\leq 0.
\end{align}

The  expression for the Mayer term in the cost~\eqref{eqn:dt_DOP_cost} is straightforward to derive, i.e.\ 
\begin{equation}
V_M^d(s_0,s_N,\pi) := V_M(\chi_0(t_0;s_0), s_N,\pi).
\end{equation}

The approximation of the integral of the running cost  $\int_{t_i}^{t_{i+1}}\ell(x(t),u(t),\theta,t) dt$  can either be computed analytically, if possible, or approximated by any suitable  numerical quadrature, i.e.
\begin{equation}
 L_i(s_i,q_i,\pi)
 :=\sum_{t\in \mathbb{T}_i^\ell} w_{i}(t) \ell\left(
\chi_i(t;s_i),\upsilon_i(t;q_i),\theta,t\right) 
\end{equation}
where the quadrature weight function $w_{i}$ is an appropriately-defined function of the finite set of points $\mathbb{T}_i^\ell\subset \mathcal{T}_i$ at which the integrand is evaluated. The quadrature scheme has to be chosen to ensure  consistency, stability and convergence of the quadrature~\cite{DR2ns,betts2010practical}; a popular choice is to evaluate the running cost at the collocation points, i.e.\ $\mathbb{T}_i^\ell=\mathbb{T}_i^f$, but other choices are possible.
\end{subequations}

The transcription is complete. The  functions in~\eqref{eqn:dt_DOP} can now be constructed from~\eqref{eqn:collocation} and the continuity constraint~\eqref{eqn:continuity}.

\subsubsection*{Existence and Uniqueness}

Even if a solution exists to the original problem~\eqref{eqn:cont_DOP}, it is possible that the discretized problem~\eqref{eqn:dt_DOP} could be infeasible, or there could be  multiple solutions to~\eqref{eqn:dt_DOP} even when the solution to~\eqref{eqn:cont_DOP} is unique. Deriving necessary or sufficient conditions for the existence and uniqueness of either problem is an on-going topic of research and is beyond the scope of this paper. However, we briefly outline here some  rules-of-thumb that often work in practice, provided the mesh is sufficiently fine. In some cases (e.g.\ for certain affine or bi-affine $f$) these rules-of-thumb are necessary and sufficient, but for  general nonlinear~$f$ and DAEs they are neither necessary nor sufficient. 

We focus only on whether a solution exists to the equality constraints~\eqref{eqn:continuity} and \eqref{eqn:collocation_f}, since this is a necessary requirement for a solution to exist and is often the main source of issues related to infeasibility and non-uniqueness.

Note that for each stage $i\in\mathbb{I}_N\backslash\{0\}$, the resulting number of equality constraints is $n_x+n_f \operatorname{card} \mathbb{T}_i^f$, where $\operatorname{card}$ denotes the cardinality of a set. 
A rule-of-thumb to avoid an over-determined set of equations is that the number of parameters for the state and algebraic variable trajectories  should be greater or equal to the number of equality constraints in~\eqref{eqn:continuity} and~\eqref{eqn:collocation_f}.  

Suppose the number of differential equations is $n_f\geq n_x+n_a\geq n_x$, where we assume that $n_a$ algebraic variables, if present, are included in the definition of~$u$, i.e.\ $n_a\leq n_u$.
Let
\[
N_i^q:= m_i/n_u-1, \ N_i^s:= n_i/n_x-1,\ N_i^f:=\operatorname{card} \mathbb{T}_i^f.
\]
Using the above rule-of-thumb we get that one should check whether
\begin{equation}
n_x(N_i^s+1)+n_a(N_i^q+1)
 \geq  n_x+n_fN_i^f
,\ \forall i \in\mathbb{I}_N\backslash\{0\}.
\label{eqn:dof}
\end{equation}

In the special case when the system is described only by an ODE, i.e.\ $n_f=n_x$ and $n_a=0$, then~\eqref{eqn:dof} reduces to 
\[
 N_i^s
\geq N_i^f,\ \forall i \in\mathbb{I}_N\backslash\{0\}.
\]
This is equivalent to saying that the number $N_i^s+1$ of basis functions   used to parameterize the state trajectory should be greater than the number of collocation points. If polynomial basis functions are used, then this implies that the degree~$N_i^s$ of the resulting state polynomials  should not be less than  the number of collocation points. 
This analysis is in agreement with the convention that in most existing polynomial collocation schemes the degree of the polynomial is chosen to be equal to the number of collocation points if the system is given by an ODE.


The analysis for the boundary equality constraints~\eqref{eqn:collocation_boundary} and stage $i=0$ proceeds in a similar manner as above, but all the equality constraints and degrees of freedom for the whole trajectory should be considered. A common rule-of-thumb to avoid an over-determined set of equations is to require that $n_E\leq 2n_x$,   but 
this is neither sufficient nor necessary, in general, especially when dealing with DAEs or a nonlinear~$f$.

In many collocation schemes the inequalities above are set to equality in order to avoid an under-determined set of equations.
Once again, this is neither necessary nor sufficient to guarantee  existence or uniqueness, in general.  It is possible to find linear and nonlinear systems for which the solutions to a collocation-based problem are non-unique even if $N_i^s
= N_i^f$. We provide two examples to demonstrate this.
\begin{example}
Suppose we wish to compute the solution to the linear initial value problem
\[
    \dot{x}(t) = \theta x(t),\ x(0) = x_0,
\]
where $\theta >0$ is given.
It possible to prove that the resulting system of linear equations from the implicit Euler, midpoint rule and trapezoidal methods all have an infinite number of solutions if, respectively, any interval length  $h_i=1/\theta$, $h_i=2/\theta$ and $h_i=2/\theta$. On the other hand, the explicit Euler method will result in a unique solution for any $h_i>0$ and~$\theta$.
\end{example}

\begin{example}
Suppose we wish to compute the solution to the nonlinear initial value problem
\[
\dot{x}(t) = \theta x(t)-x(t)^2,\ x(0) = 1.
\]
Suppose $\theta=-1$ is given. 
It possible to prove that the resulting system of nonlinear equations
from the implicit Euler, midpoint rule and trapezoidal methods all have more than one real-valued solution if, respectively, any  interval length  $h_i>0$, $h_i>0$ and $0<h_0<2$. If $h_0>2$, then there is no real-valued solution with the trapezoidal method. On the other hand, the explicit Euler method will result in a unique real-valued solution for any $h_i>0$ and~$\theta$. 
\end{example}

We therefore recommend that one proceed with caution when choosing the parameterization and collocation points, especially if $f$ is nonlinear, the dynamics are not modelled with an ODE or $n_f> n_x$. 
If the problem is infeasible, then one could increase the number of parameters used for  the state and/or free variable trajectories. If the problem is not unique, then one might wish to do the opposite. 

We note again that, in general, a collocation scheme might not be able to guarantee existence or uniqueness of a solution to~\eqref{eqn:dt_DOP}, even if a solution to~\eqref{eqn:cont_DOP} exists or is unique, respectively. However, collocation schemes have proven to be very effective in solving many challenging problems over the last few decades, so they are usually a good starting point.

\subsubsection*{Polynomial Collocation Methods}

Commonly used polynomial discretization schemes for direct collocation can be categorized into fixed-order $h$ methods and variable higher-order $p$/$hp$ methods. 
Improving the accuracy of an $h$ method is achieved by placing additional grid points during the mesh refinement process. The class of $p$/$hp$ methods, also known as pseudo-spectral methods, provide another alternative.  Improving the accuracy of a $p$ method is achieved by increasing the degree of a polynomial in one or more  intervals during the mesh refinement process; $hp$ methods allow both  an increase in  the polynomial degree and placing additional grid points. The main benefit of using  $p/hp$ methods is that, if the solution trajectories are smooth, the same accuracy can be reached with much smaller NLPs than the $h$ method counterpart, resulting in potential computational advantages. 

Well-known polynomial collocation methods\cite{betts2010practical,kelly2017,RMD2nd}  include (classified according to  degree of polynomials, collocation points and whether it is of Radau/Gauss/Lobatto type):
\begin{itemize}
    \item Explicit Euler: affine $\chi_i(\cdot;s_i)$; usually constant $\upsilon_i(\cdot;q_i)$ $\mathbb{K}_i:=\{-1\}$; Radau.
    \item Implicit Euler: affine $\chi_i(\cdot;s_i)$; usually constant $\upsilon_i(\cdot;q_i)$; $\mathbb{K}_i:=\{1\}$; Radau.
    \item Midpoint rule: affine $\chi_i(\cdot;s_i)$; usually constant $\upsilon_i(\cdot;q_i)$; $\mathbb{K}_i:=\{0\}$; Gauss.
    \item Trapezoidal: quadratic $\chi_i(\cdot;s_i)$; usually constant or affine $\upsilon_i(\cdot;q_i)$; $\mathbb{K}_i:=\{-1,1\}$,
    Lobatto.
    \item Hermite-Simpson: cubic $\chi_i(\cdot;s_i)$; usually constant; affine or quadratic $\upsilon_i(\cdot;q_i)$; $\mathbb{K}_i:=\{-1,0,1\}$; Lobatto.
\end{itemize}
In the above, we use `usually' to indicate that higher degree polynomials for $\upsilon_i$ are allowed, and might indeed be necessary, in order to guarantee that a  solution exists. This might be the case, for example,  when there are DAEs, initial and final equality constraints on the state, or other equality constraints.

If one wishes to use higher degree polynomials than cubic, care has to be taken with the choice of collocation points. This is because a uniform distribution of collocation points does not guarantee convergence as the number of mesh points $N$ increases, due to Runge's phenomenon. Instead, it might be necessary to use a non-uniform distribution of collocation points. $p/hp$ methods often use the roots of orthogonal polynomials (Legendre or Chebyshev) as collocation points.
Legendre polynomials are solutions to the Legendre differential equation defined on the interval $[-1,1]$. Within the choice of Legendre polynomials for $p$/$hp$ methods, there are three main candidates for collocation points: 
\begin{itemize}
\item Legendre-Gauss (LG) points, being the roots of a $K^\text{th}$ degree Legendre polynomial $P_{K}(\cdot)$, excludes both boundary points. 
\item Legendre-Gauss-Radau (LGR) points include $-1$ but do not include the end point at~$1$; they are the roots of the polynomial $P_{K}(\cdot)+P_{K-1}(\cdot)$. 
\item Legendre-Gauss-Lobatto (LGL) points include both boundary points $-1$ and $1$; they are the roots of $\dot{P}_{K-1}(\cdot)$. 
\end{itemize}
The most appropriate choice of  collocation points for $p/hp$ methods depends on the properties of the dynamic optimization problem. See \cite{garg2009overview,fahroo2008advances} for a more detailed discussion.

\subsection{Integrated Residual Schemes}

This is a  class of methods that  can be considered as generalization of collocation methods and is widely acknowledged to have a number of attractive properties compared to collocation methods, especially when solving high-index DAEs. These methods have been widely used for solving partial differential equations --- Galerkin methods being one of the most well-known. Here we will  outline how to generalize integrated residual methods, usually used to solve differential equations, to solving dynamic optimization problems.

The issue with collocation methods is that the differential equations~\eqref{eqn:DEs} are satisfied only at a finite set of time instances, namely the collocation points. The residual, i.e.\ the  violations of the equality constraint~\eqref{eqn:DEs}, might be non-zero elsewhere. 

As the name suggests, integrated residual methods aim to bring down the residual by formulating a constraint based on an integral of the residual over an interval. This relaxes the requirement of forcing the residual to zero at a fixed number of points, thereby allowing for the possibility of decreasing the average or maximum of the equality constraint violations.

\emph{Galerkin methods} form the basis of modern finite element methods. A generalization of Galerkin methods, namely \emph{weighted residual methods}, replaces~\eqref{eqn:collocation_f} by  the  equality constraints
\begin{equation}
    \int_\mathcal{T}
     \eta(t)^\top f(\dot{\tilde{x}}(t;s),\tilde{x}(t;s),\tilde{u}(t;q),\theta,t) dt=0, \forall \eta \in \mathcal{H},
    \label{eqn:galerkin}
\end{equation}
where each \emph{test function} $\eta:\mathcal{T}\rightarrow\mathbb{R}^{n_f}$ and $\mathcal{H}$ is a given finite set of functions. The integral in~\eqref{eqn:galerkin} can either be evaluated analytically, if possible, or approximated with any suitable quadrature. The remaining constraints and cost function are the same as with collocation.
The constraints~\eqref{eqn:galerkin} ensure that the residual is orthogonal to every test function. 
The test functions are often orthogonal to each other, but this is not necessary. 
As with collocation, care must be taken that there are enough degrees of freedom so that a solution to~\eqref{eqn:galerkin} exists. 

\emph{Least-squares methods} for initial or boundary value problems aim to solve the following problem directly or indirectly:
\begin{equation}
    \min_{q,s,\theta} \frac{1}{t_f-t_0}\int_\mathcal{T}\|f(\dot{\tilde{x}}(t;s),\tilde{x}(t;s),\tilde{u}(t;q),\theta,t)\|_2^2 dt,
    \label{eqn:ls_de}
\end{equation}
subject to the continuity constraints~\eqref{eqn:continuity} and boundary constraints~\eqref{eqn:collocation_boundary}. By comparing the first-order necessary conditions for optimality of~\eqref{eqn:ls_de} to the constraints~\eqref{eqn:galerkin},  indirect least-squares methods can be shown to be a special case of weighted residual methods, where the test functions are defined in terms of the partial derivatives of the cost function in~\eqref{eqn:ls_de}; in this case least-squares methods  are also known as Rayleigh-Ritz methods.
This observation hints at why unmodified least-squares methods, as well as other Galerkin or weighted residual methods, might lead one to conclude incorrectly that the residual and cost function in~\eqref{eqn:cont_DOP} cannot  both be  brought below given values without changing the discretization.
This is because the first-order optimality conditions of the  least-squares problem above are necessary, but not sufficient, in general. Hence, a solution to the resulting set of equations~\eqref{eqn:galerkin} could be a local maximizer of the cost in~\eqref{eqn:ls_de}.
Furthermore, these optimality conditions do not include the inequality constraints or the fact that, in a dynamic optimization problem, one is also aiming to minimise a different cost. A solution that minimizes the above least squares cost might be infeasible with respect to other inequality or equality constraints; this is also possible if  the partial derivatives with respect to $q$ and $\theta$ are not included as test functions. On the other hand, a solution that minimizes the cost in~\eqref{eqn:cont_DOP} might end up making the residual unacceptably large. 
Furthermore, the magnitudes of the residuals can only be checked \emph{a posteriori} with a weighted residual method. If the residual is too large, then the discretization mesh has to be refined, resulting in a larger optimization problem than necessary. Integrated residual methods therefore need to be modified in order to allow one to constrain or minimize the residual \emph{a priori} without having to increase the size of the optimization problem.

A  generalization of the least squares method, which allows one to solve dynamic optimization problems with inequality constraints, is as follows. The idea is to replace the equality constraints~\eqref{eqn:collocation_f} by the inequality constraints
\begin{multline}
    \int_{\mathcal{T}_i}
    \|W(t) f(\dot{\chi}_i(t;s_i),\chi_i(t;s_i),\upsilon_i(t;q_i),\theta,t)\|_2^2 dt 
    \leq h_i\varepsilon_i^f,\\
    \forall i \in\mathbb{I}_N.
    \label{eqn:least_squares}
\end{multline}
The  scalar upper bounds  $\varepsilon_i^f\geq 0$ and weight function~$W$ are assumed to be given.  
The integral in~\eqref{eqn:least_squares} can either be evaluated analytically, if possible, or approximated by any suitable quadrature. The remaining constraints and cost function are the same as with collocation. Note that the resulting $\phi_i$ is not a function of $f$, though. Note also that $\gamma_i$ is a function of both $f$ and $g$, because~\eqref{eqn:dt_DOP_ineq} is given by~\eqref{eqn:ineq_collocation} and~\eqref{eqn:least_squares}.

\begin{figure}[t]
\begin{center}
\includegraphics[width=0.8\columnwidth]{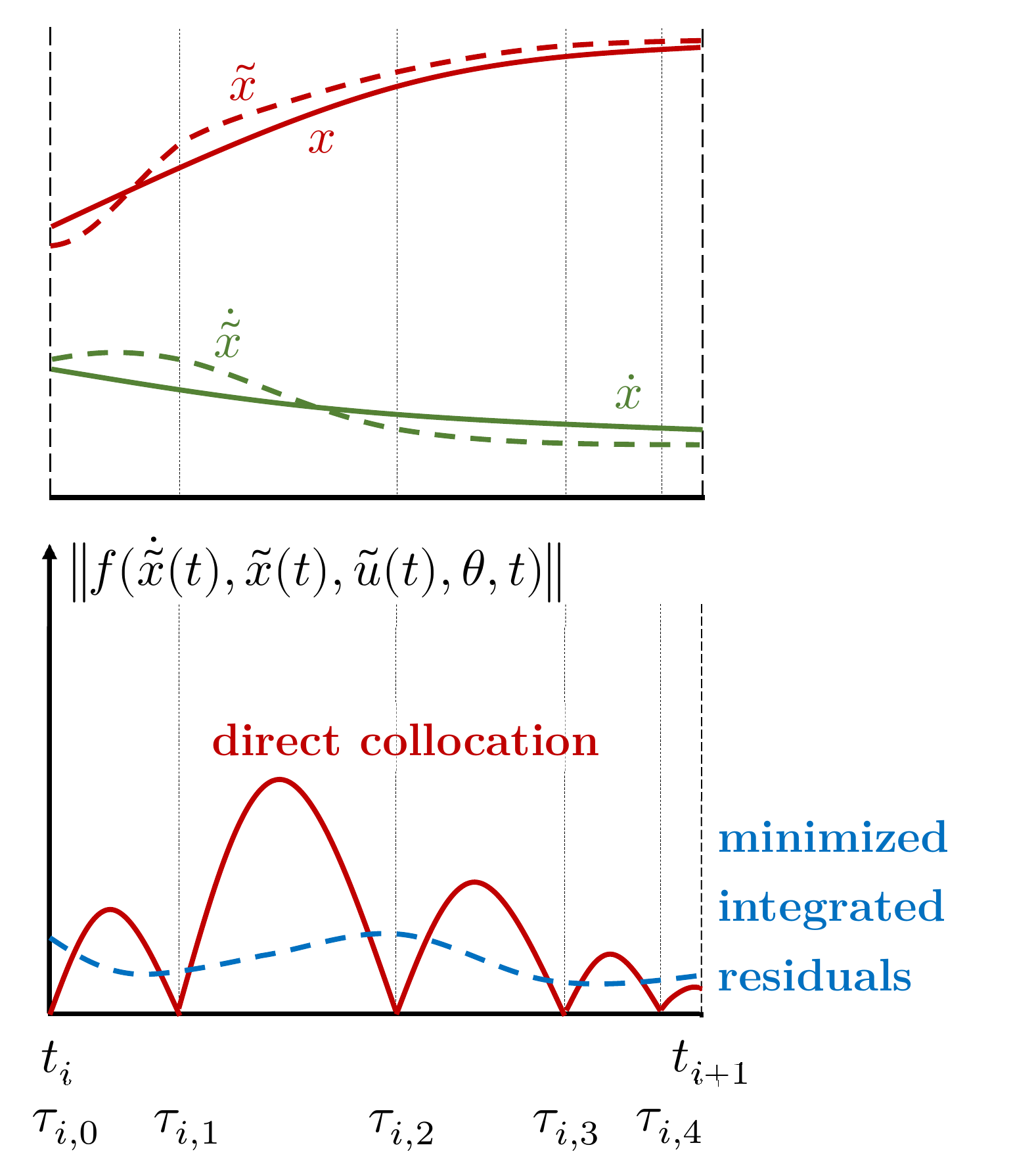}
\caption{Illustration of the differences in error characteristics between a collocation scheme and a minimized integrated residual method inside a mesh interval.} 
\label{fig:ResidualError}
\end{center}
\end{figure}

In~\eqref{eqn:least_squares},  each $\varepsilon_i^f\geq 0$ should be chosen large enough to ensure a solution to~\eqref{eqn:dt_DOP} exists. One way to compute a set of appropriate values for each $\varepsilon_i^f$ is to first solve the  weighted, nonlinear, constrained least squares problem
\[
    \min_{q,s,\pi} \frac{1}{t_f-t_0}
    \int_\mathcal{T}
    \|W(t) f(\dot{\tilde{x}}(t;s),\tilde{x}(t;s),\tilde{u}(t;q),\theta,t)\|_2^2 dt
\]
subject to~\eqref{eqn:dt_DOP_ineq}--\eqref{eqn:dt_DOP_ineq_terminal}, the continuity constraints~\eqref{eqn:continuity} and~\eqref{eqn:wp_collocation}. This solution is then used to compute bounds on the values of the integrals in~\eqref{eqn:least_squares}, before solving the resulting~\eqref{eqn:dt_DOP}. A variation of this procedure is presented in~\cite{NieKerriganCSL:2020}, 
where it is shown, via numerical examples, that this method can find more accurate solutions than collocation methods with the same mesh size and parameterization of $\tilde{u}$ and $\tilde{x}$. The difference in error characteristics between the two methods is illustrated in Figure~\ref{fig:ResidualError}.

An alternative to the above, which does not require one to first compute or estimate suitable $\varepsilon_i^f\geq 0$, is given in~\cite{NeuenhofenKerrigan:2018,NeuenhofenKerrigan:CDC2020}. However, the method in~\cite{NeuenhofenKerrigan:2018,NeuenhofenKerrigan:CDC2020} requires the availability of a suitable penalty-barrier optimization solver.

The above idea can of course be generalized and one could refer to the method as a \emph{constrained/minimized integrated residual} method, as appropriate. For example, it is not necessary to use the integral of the square of the 2-norm above --- this choice is often convenient when the NLP solver requires derivatives to exist. Suitable alternatives would be any function norm or entry-wise cost as in Section~\ref{subsec:Error_And_MeshRefinement}.


\subsection{Runge-Kutta Schemes}
\label{sec:RK}

If all the equations in~\eqref{eqn:DEs} involving derivatives of $x$ are ordinary differential equations (ODEs), i.e.\ \eqref{eqn:DEs} can be written in the semi-explicit DAE form
\begin{subequations}
\label{eqn:DAE_semi_explicit}
\begin{align}
    \dot{x}(t) &= f_{ode}(x(t),u(t),\theta,t),\ \forall t\in\mathcal{T},\label{eqn:ODE}\\
    0 &= f_{alg}(x(t),u(t),\theta,t),\ \forall t\in\mathcal{T},\label{eqn:algebraic}
\end{align}
\end{subequations}
then one can use $K$-stage Runge-Kutta schemes to perform the transcription. 

This is done by first considering the ordinary differential equation~\eqref{eqn:ODE} and forming the Runge-Kutta equations
\begin{subequations}
\label{eqn:RK_equality}
\begin{align}
    \xi_{i+1} &= \xi_i+h_i\sum_{j=1}^{K_i} b_j\tilde{f}_{i,j}  
    \label{eqn:RKx}
\end{align}
where
\begin{align}
\tilde{f}_{i,j}  
    &:= f_{ode}(\tilde{\xi}_{i,j}, \mu_{i,j},\theta,t_i+h_ic_j),
    \label{eqn:RKf}\\
    \tilde{\xi}_{i,j}&:=\xi_i+h_i\sum_{k=1}^{K_i}a_{j,k}\tilde{f}_{i,j}\label{eqn:RKxtilde}\\
    \mu_{i,j} &:= 
    \upsilon_i(t_i+h_ic_j;q_i) \label{eqn:RKu}
\end{align}
for all $i\in\mathbb{I}_N$ and $j=1,\ldots,K_i$.
The coefficients $a_{j,k},b_j,c_j$ are given by a particular scheme's Butcher tableau~\cite{betts2010practical,RMD2nd}. Recall that  $a_{j,k}=0$ for all $k\geq j$ if and only if the method is explicit, otherwise the method is implicit. It is also a good idea to define $\upsilon_i$ such that $\upsilon_i(\cdot;q_i)$ interpolates through   components of $q_i$ at the respective point $t_i+h_ic_j$, but this is not necessary. 

All polynomial collocation schemes  can be shown to be equivalent to a Runge-Kutta scheme, where the set of collocation points are given by $\mathbb{T}_i^f=\{t_i+h_ic_j \mid j=1,\ldots,K_i\}$, with all $c_j$ distinct~\cite{betts2010practical,kelly2017,RMD2nd}, i.e.\ $\mathbb{T}_i^f$ has $K_i$ elements. However, not all Runge-Kutta methods are equivalent to polynomial collocation schemes, e.g.\ when all the~$c_j$ are not distinct.
Hence, it may not be possible to find an interpolating polynomial such that $\chi_i$ interpolates through all estimates of the state $\tilde{\xi}_{i,j}$. We therefore assume here that the chosen Runge-Kutta scheme is \emph{not} equivalent to a polynomial collocation method. Care must therefore be taken when choosing and interpreting appropriate parameterizations for certain classes of Runge-Kutta schemes.

A suitable choice for $\chi_i$ and $s_i$ is to satisfy
\[
    \chi_i(t_i;s_i)= \xi_i \text{ and }   \chi_i(t_{i+1};s_i)= \xi_{i+1}, \ \forall i \in\mathbb{I}_N.
\]
One could then choose $s_i:=(\xi_i,\xi_{i+1})$ and  $\chi_i(\cdot;s_i)$ to be affine.
However, this is not always desirable if one wishes to use sparsity-exploiting NLP solvers to compute a solution~\cite{betts2010practical}. If a sparse NLP solver is available, then one could choose to also include some or all of the estimates of the state $\tilde{\xi}_{i,j}$, state derivative $\tilde{f}_{i,j}$ and free variables  $\mu_{i,j}$ in the definition of $s_i$ and add the relevant parts of~\eqref{eqn:RKf}--\eqref{eqn:RKu} as equality constraints to~\eqref{eqn:dt_DOP_eq}. One also needs to define~$\chi_i$ appropriately, but with the above-mentioned limitation on  being able to interpolate only some, but not all, of the state estimates.

The algebraic constraints~\eqref{eqn:algebraic} can be enforced by adding the following constraints to the above:
\begin{equation}
  0 = f_{alg}(\chi_i(t;s_i), \upsilon_i(t;q_i),\theta,t),\ \forall t\in\mathbb{T}_i^a,
\end{equation}
where $\mathbb{T}_i^a$ is any finite subset of $\mathcal{T}_i$. 
\end{subequations}
The equality constraints~\eqref{eqn:dt_DOP_eq} are then given by~\eqref{eqn:RK_equality}, the continuity constraints~\eqref{eqn:continuity} and~\eqref{eqn:wp_collocation}.

The inequality  constraints~\eqref{eqn:ineqs}
can be enforced in a similar manner as for collocation methods, i.e.\ \eqref{eqn:dt_DOP_ineq} is given by \eqref{eqn:ineq_collocation}. However, for constraints that involve the derivative of the state, i.e.\ $\dot{\chi}_i$, one could use  $\tilde{f}_{i,j}$ or $f_{ode}(\chi_i(t;s_i),\upsilon_i(t;q_i),\theta,t)$, depending on the choice of Runge-Kutta scheme and definition fors $\chi_i$ and $s_i$. 

The remaining constraints and cost function terms in~\eqref{eqn:dt_DOP} can also be derived in a similar manner to collocation. 
The most straightforward choice for the points at which constraints are enforced or the running cost is evaluated is to choose
$\mathbb{T}_i^a=\mathbb{T}_i^g=\mathbb{T}_i^\ell=\mathbb{T}^c\cap\mathcal{T}_i=\{t_i,t_{i+1}\}$, but of course other finite subsets of $\mathcal{T}_i$ are possible, provided care is taken as above.
In particular, one has to ensure that the parameterization has been chosen such that there are enough degrees of freedom for a solution to exist, which is not always straightforward for systems described by DAEs~\cite{betts2010practical,kunkelMehrmann2006}. 

\subsection{Shooting Methods} \label{sec:shooting}
We briefly discuss another very popular class of methods for solving dynamic optimization problems, namely \emph{shooting} schemes. The word `shooting' describes the process whereby a solution at a later time-step is integrated from available solutions at one or more previous time steps, a procedure also known as time-marching. 
Dynamic optimization methods based on time-marching are commonly known as \emph{sequential methods}, where initial states, parameters and free variables are iteratively adjusted,  with the help of sensitivity information, until all path constraints and boundary conditions are satisfied. For more details, see~\cite{betts2010practical,RMD2nd,iet:/content/journals/10.1049/iet-cta.2019.0168}.

When this approach is implemented using integration from an initial state all the way to the final time, the method is known as \emph{single shooting}. Practical use of single shooting methods generally requires the dynamics to be stable. Furthermore, the method can be very sensitive to numerical inaccuracies and initial guesses, leading to unstable and ill-conditioned boundary value problems (BVPs). Consequently, the solution process tends to have a high chance of failure. 

One way to mitigate the shortcomings of single shooting is to subdivide the grid into multiple intervals that are connected with corresponding continuity conditions. As a result, time marching only needs to be implemented on a short time interval. Doing so allows for the application of shooting methods to unstable systems and makes the method much more robust to numerical inaccuracies. This approach is known as \emph{multiple shooting}, and it could be considered to be a hybrid between a simultaneous and sequential method, because the state and input trajectories must be solved altogether as a whole to yield a valid solution.

Multiple shooting ensures that the differential equations~\eqref{eqn:DEs} are approximately satisfied (up to a specified tolerance) over  an interval~$\mathcal{T}_i$ as follows. Suppose the state $\tilde{x}(t_i):=s_i$ at time $t=t_i$ and a trajectory of free variables $\upsilon_i(\cdot;q_i)$ is given. Any suitable, stand-alone differential equation solver (including adaptive variable-step/variable-order solvers) is used to compute $\tilde{x}(t_{i+1}):=\varphi_i(t_{i+1};s_i,q_i,\theta,t_i)$, which is the evaluation at $t=t_{i+1}$ of an approximate solution to the differential equations~\eqref{eqn:DEs}, with $u(\cdot):=\upsilon_i(\cdot;q_i)$ and initial condition $x(t_i)=s_i$. The continuity constraints~\eqref{eqn:continuity} are then replaced by
\[
    s_{i+1} = \varphi_i(t_{i+1};s_i,q_i,\theta,t_i),\ \forall i\in\mathbb{I}_N
\]
with no other constraints being functions of~$f$.
Note that the analytical expression for $\varphi_i$ is not actually computed. Instead, $\varphi_i(t_{i+1};s_i,q_i,\theta,t_i)$ is the output of an algorithm with  $(f,\upsilon_i,s_i,q_i,\theta,t_i,t_{i+1})$ as its input.

The remaining constraints and cost function terms are then obtained in a similar fashion as above for collocation or  Runge-Kutta methods, where  $\chi_i(t_i;s_i)$ and $\chi_i(t_{i+1};s_i)$ are replaced by $s_i$ and $\varphi_i(t_{i+1};s_i,q_i,\theta,t_i)$, respectively, for all $i\in\mathbb{I}_N$. Note that if~\eqref{eqn:DEs} is given by implicit DAEs and there are inequality constraints on some of the  state derivatives, then estimates of the state derivatives should also be provided by the DAE solver; if the system is a semi-explicit DAE or ODE then one could use~\eqref{eqn:ODE} to estimate state derivatives.

In order to simplify implementation, often the mesh and other time instances are chosen such that $\mathbb{T}_i^g=\mathbb{T}_i^\ell=\mathbb{T}^c\cap\mathcal{T}_i=\{t_i,t_{i+1}\}$, but other choices for finite subsets of $\mathcal{T}_i$ are possible. For example, in many cases the  interface to the differential equation solver can provide approximate solutions $\varphi_i(t;s_i,q_i,\theta,t_i)$ for all~$t$ in a given finite subset of~$\mathcal{T}_i$. In this case, one can see that single shooting is a special case of multiple shooting with $N:=1$, where $\mathbb{T}_0^g\cup\mathbb{T}_0^\ell\cup\mathbb{T}^c$ is a given finite subset of $\mathcal{T}$ and $t\mapsto\tilde{u}(t;q)$ is a discontinuous piecewise differentiable trajectory, parameterized by $q=q_0$.

\subsection{Solving the NLP}
\label{sec:NLP}

The choice of the appropriate solver for the transcribed problem~\eqref{eqn:dt_DOP} depends on the characteristics of the optimal control problem and available computational resources. It is beyond the scope of this paper to go into this fascinating topic, which is a highly active area of research, hence we only briefly highlight some of the main points to consider.

One may choose between NLP solvers that are derivative-free, solvers that use first-order derivative information, solvers that use first- and second-order information, or solvers based on first derivatives and a quasi-Newton approximation of the second derivative information (e.g.\ BFGS algorithms). 
Generally speaking, when the solution is smooth and well-behaved, solvers which use second-order derivative information will converge in fewer iterations than first-order or derivative-free methods. However, this derivative information may need to be obtained through sparse finite differences, analytical derivatives or algorithmic differentiation packages. These can be tedious to derive or  expensive to compute. Conversely, solvers that only use first-order or no derivative information may need a larger number of iterations to converge, but the computational effort per iteration is typically significantly lower. In practice, there is no class of NLP algorithm that is best suited for all problem types. A proper choice should be made on a case to case basis. 

When using derivative-based solvers, fast and accurate computation of derivative information is key to solving the problem efficiently. Different orderings of optimization decision variables lead to different sparsity patterns of the Jacobian and the Hessian of the relevant functions. Exploiting the sparsity patterns of these matrix systems can lead to substantial reductions in computational complexity and storage usage, in both the process of supplying the derivative information and solving the NLP. In fact, sparse linear algebra has become one of the most important aspects of numerical optimal control, allowing the efficient solution of large-scale practical problems. 

Although an NLP solver tailor-made for specific transcription methods can appear attractive~\cite{FERREAU201713194,RMD2nd}, in many cases off-the-shelf NLP solvers can be used directly. Popular candidates include IPOPT~\cite{IPOPT} (an interior point solver), SNOPT \cite{SNOPT} (a  sequential quadratic programming (SQP) solver using an active-set quadratic programming (QP) solver), WORHP \cite{WORHP} (an SQP solver using an interior point QP solver), and NOMAD \cite{NOMAD} (a derivative-free solver).

\begin{figure}[tb]
\begin{center}
\includegraphics[width=\columnwidth]{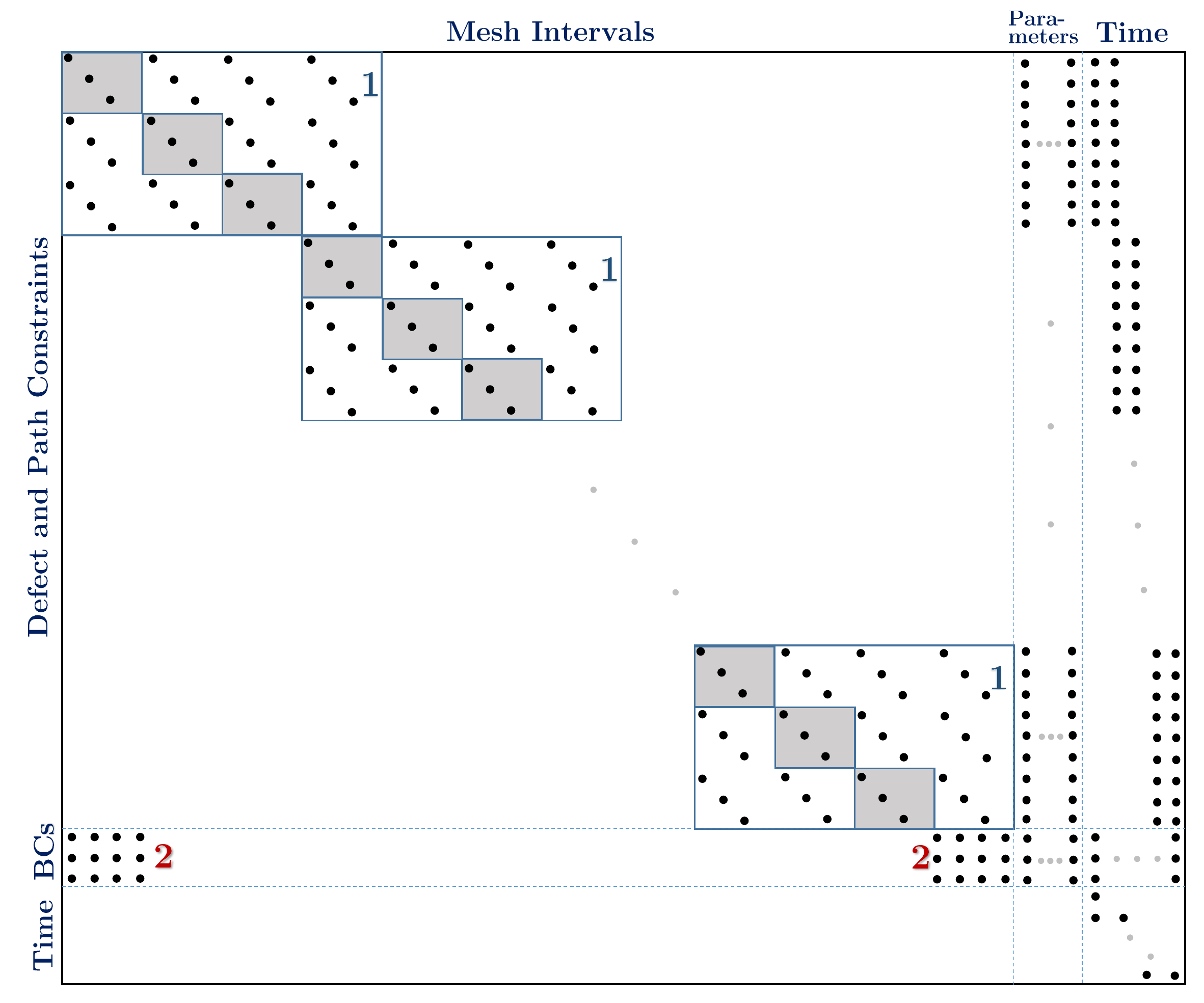}
\caption{Sparsity pattern for the constraint Jacobian using LGR collocation when grouped by stage (Block Type 1: block structures also with internal block structures; Block Type 2: from boundary conditions, either zeros or endpoints)} 
\label{fig:LGRReorderJacobianSparsity}
\end{center}
\end{figure}

\begin{figure}[tb]
\begin{center}
\includegraphics[width=\columnwidth]{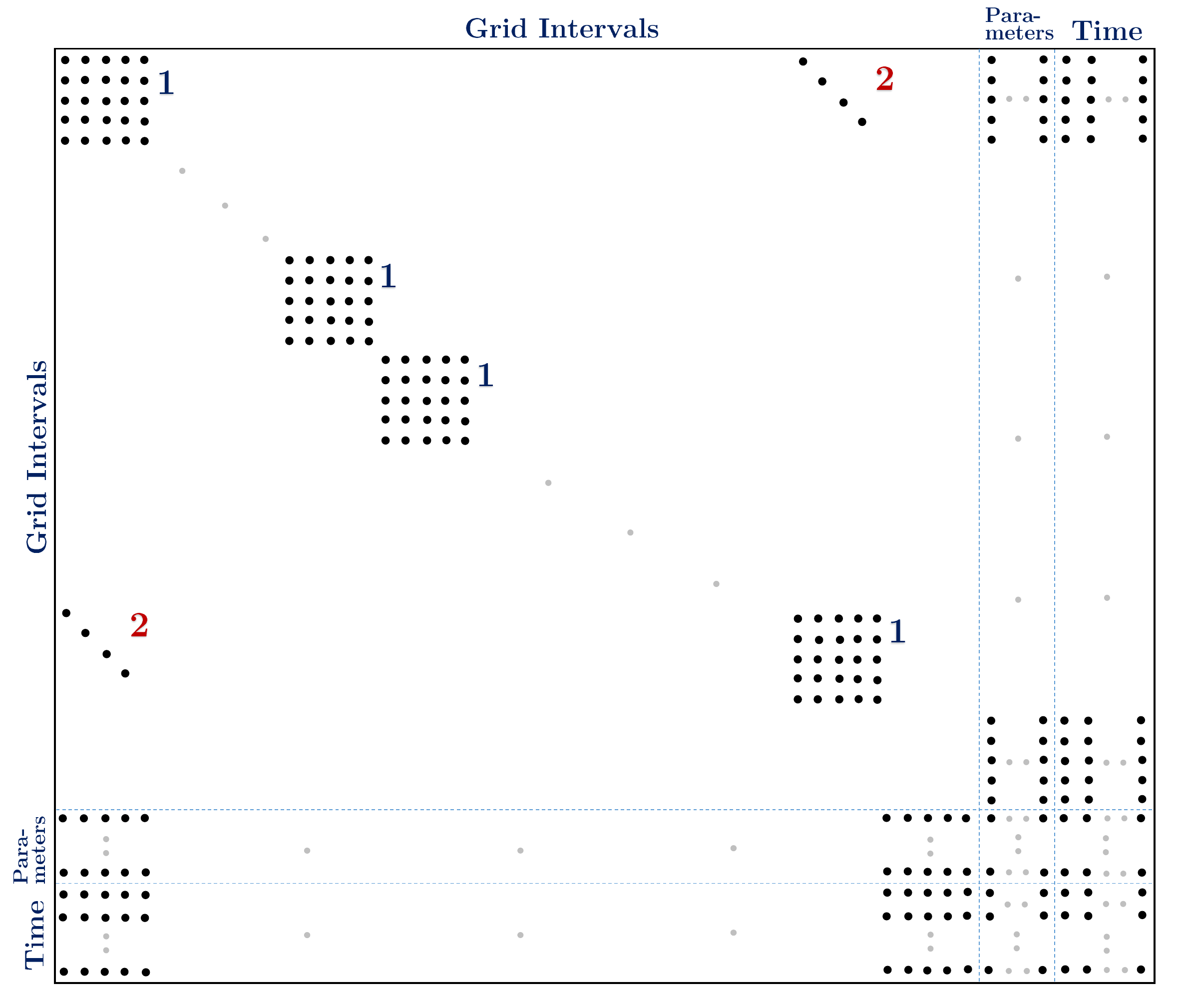}
\caption{Sparsity pattern for the Hessian of the cost function using LGR collocation when grouped by stage (Block Type 1: populated block; Block Type 2: diagonal block)} 
\label{fig:LGRReorderHessianSparsity}
\end{center}
\end{figure}

We briefly outline here the main reason why the structure of~\eqref{eqn:dt_DOP} gives rise to sparse matrices with exploitable structure in derivative-based NLP solvers. Suppose we order the sequence of decision variables as 
$
    (s_0,q_0,s_1,q_1,\ldots,s_{N-1},q_{N-1},s_N,\pi)
$
where the $s_i$ and~$q_i$ are ordered by increasing stage number. 
Suppose also that the constraints are ordered in increasing stage number, with the boundary constraints after the last stage. It can then be shown that the rows of the Jacobian of the constraints
can be permuted to be a bordered block-banded  matrix (also known as a block arrowhead matrix)~\cite{betts2010practical}, as  in Figure~\ref{fig:LGRReorderJacobianSparsity}.  It can also be shown that the Hessian of the cost function is a symmetric, bordered block-diagonal matrix, as  in Figure~\ref{fig:LGRReorderHessianSparsity}. The above two facts can be used to show that the KKT matrix  can be permuted to be a symmetric, bordered block-banded matrix. Efficient linear algebra solvers exist that can exploit this structure, while guaranteeing that the computational complexity scales linearly with the number of stages $N$~\cite{NWKSB}. 


An alternative is to define the new decision variable $z_0:=\pi$, introduce  new decision variables $z_i$ and  add the
constraints $z_{i+1}=z_i$ for all $i\in\mathbb{I}_N$. The decision variables can then be defined and ordered as
$(s_0,z_0,q_0,\ldots,s_{N-1},z_{N-1},q_{N-1},s_N,z_N)$. The KKT matrix of this lifted problem can be permuted to be block-banded. This opens up a potentially larger number of efficient structure-exploiting solvers~\cite{RMD2nd,FERREAU201713194,NWKSB}, but at the expense of having a bigger optimization problem.






\section{Case Study: Multi-patient Ventilation}
\label{sec:ventilators}

COVID-19 is a viral illness caused by a newly discovered coronavirus, named Severe Acute Respiratory Syndrome Coronavirus~2 (SARS-CoV-2). The disease was first documented in Wuhan, Hubei Province, China with a number of unexplained pneumonia cases in December 2019~\cite{doi:10.1002/jmv.25722}. The disease has subsequently rapidly spread worldwide, leading to drastic measures to be taken globally to try to limit further spread and contain the infection. 
On 18 September 2020, there had been 30,055,9710 confirmed cases and 943,433 confirmed deaths from COVID-19 reported worldwide, with 216 countries, areas or territories affected~\cite{WHO}. The virus continues to spread with a reproductive number estimated by the WHO as 2.5 (higher than for influenza).

COVID-19 causes a multitude of symptoms with the main ones being fever and a dry cough. Although many cases are asymptomatic, severe illness can cause death, the risk of which increases with age and certain pre-existing co-morbidities. The main cause of death is respiratory failure, with myocardial damage and circulatory failure also contributing~\cite{CCDCP}. 
In the prevention and emergency handling of respiratory failure, patients may require positive pressure mechanical ventilation and in severe cases, the use of extracorporeal membrane oxygenation (ECMO) treatment.

As a result, the availability of ventilators and ECMO machines may become a decisive factor in outcomes for many patients with severe disease. At the time of writing, there was worldwide concern that there will be a shortfall of intensive care beds and of mechanical ventilators to support the most severe cases of COVID-19. Estimates of the number of ventilators in the US in March 2020 ranged from 60,000 to 160,000~\cite{JHCHS}, although the distribution of ventilators is unlikely to directly coincide with COVID-19 hotspots. There were varying estimates for the number of ventilators that could be required by the US at the height of the pandemic (potentially up to 1~million), but whichever estimate is used, there was concern that the national strategic reserve would be insufficient to fill the projected gap~\cite{doi:10.1056/NEJMp2006141}.

These concerns have lead to suggestions of how to bridge the gap, including ventilating multiple patients using one ventilator. This method has been widely debated and various problems have been identified, including the risk of infection, the inability to deliver different pressures or achieve different tidal volumes in individual patients (with volume being delivered to the most compliant lungs) and difficulties with patient monitoring. The American Society of Anaesthesiologists released a statement in March 2020 advising against ventilating multiple patients per ventilator (while any clinically proven, safe and reliable therapy remains available). They give multiple reasons for this, including the above, as well as that positive end-expiratory pressure (PEEP), which is of critical importance to these patients, would be impossible to manage, difficulties with monitoring of pulmonary mechanics and alarm monitoring, difficulties arising from one patient deteriorating suddenly or having a cardiac arrest and ethical issues~\cite{ASA}.

In effect, the practice of splitting a single ventilator to service multiple patients is seen as experimental and untested, with little to no sound principles/methods on how to safely manage each patient. Preceding the COVID-19 crisis, there had been only few numerical or experimental studies on the effect of split ventilation. A brief summary of pre-existing experimental studies may be found in \cite[Table~1]{tronstad2020splitting} and consist primarily of tests conducted on mechanical lungs and animals, although anecdotal reports do tell of successful split ventilation of humans in past crises. Further research is therefore needed into the feasibility of ventilating multiple patients using one ventilator if this is to become a viable option. Given the current global crisis, or potential of any other future similar pandemics, the solution to this problem could be of vital importance.

\subsection{Modeling of the Ventilator-Patient System}
Employing the analogue between electrical current and air flow rate, it is common to model human lungs as series connections of resistors and capacitors~\cite{campbell1963electrical}. The capacitor represents the combined compliance of the lung and thorax, which may vary significantly even between normal adults.  Resistors model pressure losses due to the restriction of airflow through the ventilation inhale/exhale pipes and endotracheal tube. Finally the pressures of the ventilators are modelled by, potentially time-varying, voltage sources.

The control inputs, namely those variables that can be directly manipulated/actuated, in single patient ventilation are the peak inspiratory pressure (PIP) and PEEP pressures generated at the ventilator,  respiratory rate, and the inhale ratio --– all of which must be manipulated within medically safe ranges to ensure the minimum and maximum pressure and the total air volume received by the patient are appropriate for the treatment. A key issue in delivering clinician-prescribed tidal volumes to multiple patients is that the division of air flow between patients is largely insensitive to these normal control inputs. Instead, the ratio of flow going to each patient is highly dependent on the resistance and lung compliance of the individual patients. To this end, as proposed in~\cite{doi:10.1098/rsos.200585}, we add variable resistances to each patient circuit, which may be used to change the effective resistance of each patient as seen by the ventilator.

\begin{figure}[t]
\begin{center}
\includegraphics[width=\columnwidth]{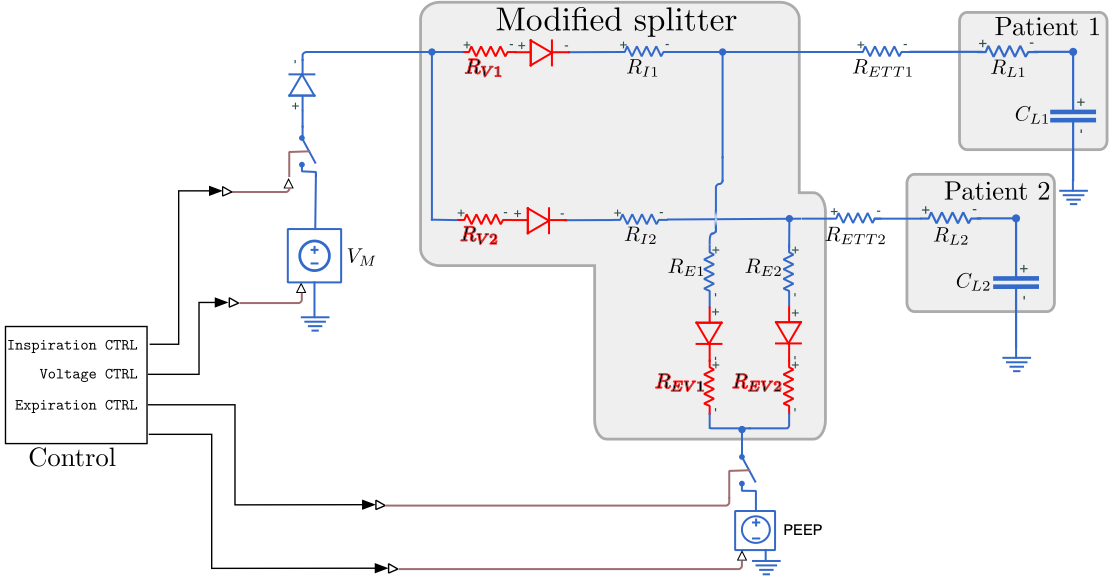}
\caption{Diagram of electrical circuit analogue of the splitter+patient system. The adjustable resistances in the inhale and exhale paths, as well as added  check valves, are indicated in red with subscripts $\cdot_{Vp}$ and diodes, respectively. The exhale pressure is controlled by the ventilator via the block labeled $V_M$, and  the exhale pressure via the block labeled PEEP. $R_{Ip}$ and $R_{Ep}$ are the resistances of the tubes in the inhale and exhale paths, respectively, and $R_{ETTp}$ is the resistance of the endotracheal tube.} 
\label{fig:CircuitModel}
\end{center}
\end{figure}

This circuit model is formalized in Figure~\ref{fig:CircuitModel} for two-patients. 
The adjustable resistances, shown in red, are added to both inhale and exhale paths for each patient and are assumed known and manipulable. It is not possible to distinguish between unknown series resistances in the estimation problem --- and indeed doing so offers no obvious benefit to control. The unknown series resistances in the inhale and exhale path of a patient $p$ are denoted as $R_p^I := R_{Ip} + R_{ETTp} + R_{Lp}$ and $ R_p^E := R_{Ep} + R_{ETTp} + R_{Lp}$.

To capture higher-order nonlinear dynamics, which may not be accurately described by linear circuitry, we will introduce fictitious `quadratic' resistances into the estimation and control problems below. This is motivated by  experimental results for a small orifice restrictor~\cite{Plummer2020.04.12.20062497}, which show that a quadratic fit to mean flow rate vs pressure drop data results in smaller errors than a linear fit.

The lungs of a patient with acute respiratory distress syndrome can have time-varying resistances and compliances. In particular, the resistance between the inhale and exhale phases could be different~\cite{doi:10.1164/rccm.201612-2495CI}. To keep the discussion simple in this tutorial, we will assume that the patient parameters are constant. However, extensions to model time-varying patient parameters are possible within the general dynamic optimization framework presented in this paper.

Notably, this is a human-in-the-loop system, where both medical practitioners and patients form part of the closed-loop system. In short, ventilated patients are housed in intensive care units, where clinicians  monitor their vital signs, including heart rate and oxygen saturation. Based on these measurements the clinician may change the ventilator mode of operation, the reference tidal volume, the value of the adjustable resistances and/or other set-points.

\subsection{System Identification and Estimation} \label{sec:CaseStudy_Estimation}

Suppose there are $n_p$ patients connected to a ventilator. Inhalation occurs for $t_I$ seconds during the time interval $T_I:=[t_0,t_0+t_I)$, followed with exhalation for $t_E$ seconds during the time interval $T_E:=[t_0+t_I,t_f)$, where $t_f:=t_0+t_I+t_E$.  The inhale to exhale ratio is therefore $\rho:=t_I/t_E$ and the respiratory rate is  $f_B:=60/(t_I+t_E)$ breaths per minute.

The PIP~$V_I(\cdot)$ and PEEP~$V_E(\cdot)$ are known. We also assume that the linear coefficient vector $R_\delta$ and quadratic coefficient vector $R_\delta^Q$  of  adjustable resistances are known. The adjustable resistances can be manipulated between zero and a maximum value, which we denote as fractions by $a_{I,p}\in[0,1]$ and $a_{E,p}\in[0,1]$ for the inhalation and exhalation tubes, respectively, for patient $p$.
Depending on what is of interest to the clinician, we would, for example, aim to obtain estimates or upper and lower bounds for the volume of air delivered to the patient, resistance values in the set $\mathcal{R}:=
\cup_{p\in\mathbb{P}}\{R^I_{p},R^E_{p},R^{IQ}_p,R^{EQ}_p\}$ 
and/or compliance values in the set $\mathcal{C}:=\cup_{p\in\mathbb{P}}\{C_p\}$, where $\mathbb{P}:=\{1,\ldots,n_p\}$. 


Suppose we have a sequence $\mathcal{Y}$ of vectors of noisy measurements at $n_m$ distinct time instances $t_k\in T_B$, $k\in \mathbb{K}:=\{1,\ldots,n_m\}$, during the time interval $T_B:=[t_0,t_f)$ of one breath, i.e.
\[
    \mathcal{Y}:=\left(y(t_1),\ldots,y(t_{n_m})\right), 
\]
where each $y(t_k)$  consist 
of the flow rate \emph{out} of the ventilator $y_0(t_k)$ (if the flow rate is negative, then flow is into the ventilator), i.e.
\begin{subequations}
\label{eqn:measure_flow_rate}
\begin{align}
y_0(t_k) = \nu_{0,k}+ \sum_{p\in \mathbb{P}} i_p(t_k),
\ 
    |\nu_{0,k}| \leq \overline{\nu}_0,\ 
    \forall k \in \mathbb{K},
    \label{eqn:measure_flow_rate_ventilator}
\end{align}
as well as measurements of the individual  flow rate \emph{into} a patient $y_p(t_k)$ (if the flow rate is negative, then flow is out of the patient) for some subset of patients $\mathbb{P}_m\subseteq \mathbb{P}$, i.e.
\begin{align}
y_p(t_k) = \nu_{p,k} +  i_p(t_k),
\ 
    |\nu_{p,k}| \leq \overline{\nu}_p,\ 
    \forall k \in \mathbb{K},\ p \in\mathbb{P}_m
    \label{eqn:measure_flow_rate_patient},
\end{align}
\end{subequations}
where $\nu_{\cdot,k}$ denotes measurement noise at time $t=t_k$.
In other words, each $y(t_k)$ is composed of $1+\operatorname{card} \mathbb{P}_m$ measurements, appropriately ordered. The vector of measurement noises $\nu$ and the vector of noise bounds $\overline{\nu}$ are defined in a similar manner to $\mathcal{Y}$.


We are interested in characterizing the set of all possible parameters $\mathcal{R},\mathcal{C}$, pressures $v$, flow rates $i$ and disturbances~$w$ that are consistent with noisy measurements~\eqref{eqn:measure_flow_rate} and that satisfy the following  equations $\forall p\in \mathbb{P}$:
\begin{subequations}
\label{eqn:dynamics}
\begin{align}
    (C_p, R^I_{p}, R^E_{p},R^{IQ}_p, R^{EQ}_p) \geq 0, \ &\\
    i_p(t)  = C_p\dot{v}_p(t),\ 
    v_p(t) \geq 0,\ |w_p(t)| \leq \overline{w},\ \forall t &\in T_B\\
    \delta_p(t) := R_\delta i_p(t)+\operatorname{sign}(i(t))R_\delta^Q i_p(t)^2,\ \forall t &\in T_B\\
    R^I_{p}i_p(t)+R^{IQ}_pi_p(t)^2+a_{I,p}\delta_p(t) =& \notag \\
    V_I(t)-v_p(t)+w_p(t)\ \forall t &\in T_I\\
    i_p(t) \geq 0,\ \forall t &\in T_I\\
     R^E_{p}i_p(t)-R^{EQ}_pi_p(t)^2 +a_{E,p}\delta_p(t)= & \notag \\
    V_E(t)-v_p(t)+w_p(t),\ \forall t &\in T_E\\
    i_p(t) \leq 0,\ \forall t &\in T_E
\end{align}
where the pressures in the lung $v:T_B\rightarrow \mathbb{R}^{n_p}$ are  continuous, but the flow rates $i:T_B\rightarrow \mathbb{R}^{n_p}$ and disturbance signals $w:T_B\rightarrow \mathbb{R}^{n_p}$ can be discontinuous. In the above,  the measurement noise $\nu:=(\nu_1,\ldots,\nu_{n_m})\in\mathbb{R}^{n_m}$ is assumed to be bounded by~$\overline{\nu}$.
The disturbance signals $w(\cdot)$  are bounded by 
$\overline{w}$ and represent modelling errors due to 
neglecting some dynamics. 
Of course, if there is no measurement noise or modelling error, we could set   $\overline{\nu}$ or $\overline{w}$, respectively, to zero  or just remove $\nu$ or $w$  from the unknowns. 
Note that one could also use substitution to eliminate $\nu$ and $w(\cdot)$, which would result in a smaller optimization problem. We could also eliminate $i(\cdot)$ if there are no quadratic terms in $i(\cdot)$.

We will assume that the patients are in (quasi) steady-state in the sense that the trajectories are periodic. We can then add the $n_p$ constraints
\begin{align}
    v(t_0) = v(t_f). \label{eqn:limit_cycle}
\end{align}
\end{subequations}

Let $\mathcal{Z}(\mathcal{U})$ be the set of  parameters $\vartheta:=(\mathcal{R},\mathcal{C})$ and  functions $\sigma:=(v,i,w)$ that satisfy~\eqref{eqn:dynamics} for the given manipulated variables
\[
\mathcal{U}:=(V_E,V_I,a_I,a_E,t_I,t_E).
\] 
To be precise,
\[
\mathcal{Z}(\mathcal{U}):=\{z:=(\vartheta,\sigma) \mid \text{\eqref{eqn:dynamics} is satisfied } \forall p \in\mathbb{P}\}.
\]
The set of  unknowns that are consistent with the measurements is
\[
\Omega(\mathcal{U},\mathcal{Y}) := \{w := (z,\nu) \mid z \in \mathcal{Z}(\mathcal{U})\ \text{and \eqref{eqn:measure_flow_rate} 
is satisfied}\}.
\]

To provide an estimate of the unknowns $\omega$, a suitable optimization problem to solve, which is  widely used in the estimation and system identification literature~\cite{RMD2nd,betts2010practical}, is
\begin{align}
    \hat{\omega}(\mathcal{U},\mathcal{Y}):=\underset{\omega\in\Omega(\mathcal{U},\mathcal{Y})}{\operatorname{argmin}}\ \nu^\top S_\nu^{-1}\nu +  
    \int_{t_0}^{t_f}w(t)^\top S_w^{-1}w(t) dt \label{eqn:ls}
\end{align}
where 
$S_\nu,S_w$ are given positive definite matrices,
so that the estimates of the unknowns can be found from $$\hat{\omega}=:(\hat{\vartheta},\hat{\sigma},\hat{\nu}).$$ This problem can be interpreted as a weighted, constrained least squares fit of the measurements to the differential and algebraic equations, with $S_\nu,S_w$ the tuning variables used to trade-off the error due to  measurement noise and neglected dynamics.

One of the key variables that the clinicians  care about is the amount of air inhaled and exhaled by the patient. 
In other words, we need to provide  estimates and/or bounds for the \emph{tidal volume}, defined as the amount of air exhaled by patient~$p$, i.e.\
\begin{multline*}
    \Delta Q_E(i_p):= 
    -\int_{t_0+t_I}^{t_f}i_p(t)dt
    =C_p(v_p(t_0+t_I)-v_p(t_f))
\end{multline*}
for each $p\in\mathbb{P}$. 
Note that the amount of air exhaled is the same as the amount of air inhaled when the patient is at steady-state, i.e. $ \Delta Q_E(i_p)=\Delta Q_I(i_p) := \int_{t_0}^{t_0+t_I}i_p(t)dt
$. 


\subsection{Control}
Suppose we are given some desired tidal volume $\Delta Q_p^d$ to be delivered to patient $p$ with maximum tolerance $\epsilon_p>0$.  
We first proceed with using the estimated parameters from solving~\eqref{eqn:ls}, under the assumption that there are modelling errors. However, to simplify the presentation in this paper,  we adopt the commonly-adopted practice of assuming certainty equivalence in the control problem, i.e.\ $w = 0$; the extension to a robust/stochastic control formulation is beyond the scope of this tutorial paper.

Given estimates of the unknown parameters $\hat{\vartheta}(\mathcal{U}_{old},\mathcal{Y})$ obtained with the  measurements $\mathcal{Y}$ and manipulated variables~$\mathcal{U}_{old}$, we  proceed to compute new manipulated variables~$\mathcal{U}_{new}$ by solving the  optimal control problem: 
\begin{subequations} \label{eqn:Control_Problem}
\begin{align}
    (\mathcal{U}_{new}^*(\mathcal{U}_{old},\mathcal{Y}),z^*(\mathcal{U}_{old},\mathcal{Y}))\in \underset{(\mathcal{U}_{new},z)}{\operatorname{argmin}}\  J(\mathcal{U}_{new},z)
\end{align}
subject to
\begin{align}
    z \in \mathcal{Z}(\mathcal{U}_{new}),\ \vartheta \in \hat{\vartheta}(\mathcal{U}_{old},\mathcal{Y}),\ &\\
    |\Delta Q_E(i_p) - \Delta Q_p^d|\leq \epsilon_p,\ &\forall p \in \mathbb{P},\\
    \underline{V_I} \leq V_I(t) \leq \overline{V_I},\ &\forall t\in T_B \\
    \underline{V_E} \leq V_E(t) \leq \overline{V_E},\ &\forall t\in T_B\\
    w(t) = 0, \  & \forall t \in T_B\\
    0 \leq a_{I,p}\leq 1,\ & \forall p \in \mathbb{P}\\
    0 \leq a_{E,p} \leq 1,\ & \forall p \in \mathbb{P}\\
    \underline{f_R} \leq f_R \leq  \overline{f_R},\ &\\
    \underline{\rho} \leq \rho \leq  \overline{\rho},\ &\\
    0 \leq t_I,\ 0 \leq t_E.\ &
\end{align}
\end{subequations}
In the above $\overline{\cdot}$ and $\underline{\cdot}$ are given upper and lower bounds on the respective variables.
We consider minimising the energy used by the ventilator, i.e.
\[
J(\mathcal{U}_{new},z):=\int_{ T_I}V_I(t)i(t)dt+\int_{ T_E}V_E(t)i(t)dt
\]
where the rate of air flow out of the ventilator is 
$
i(t):=\sum_{p\in\mathbb{P}} i_p(t)$.

\section{Simulation Results}
\label{sec:results}
The problems were transcribed using the  dynamic optimization toolbox \texttt{ICLOCS2} \cite{nie2018iclocs2} in MATLAB. The NLP was solved with interior point solver \texttt{IPOPT} \cite{IPOPT} to a relative convergence tolerance (\texttt{tol}) of $10^{-9}$. In \texttt{ICLOCS2}, both state and free variable trajectories are continuous trajectories inside a single phase, but allowed to be discontinuous with a multi-phase setup.

\subsection{System Identification and Estimation}
\begin{figure}[t]
\begin{center}
\includegraphics[width=\columnwidth]{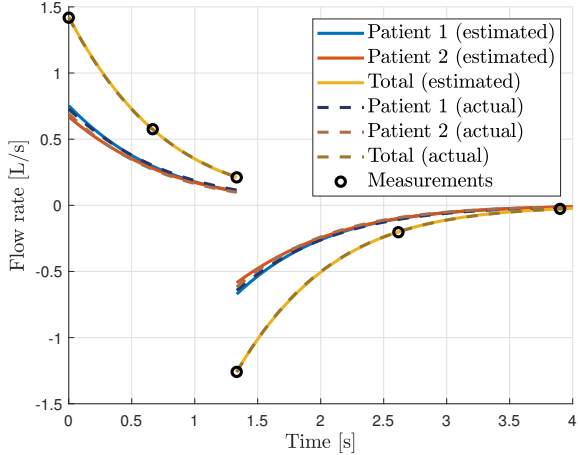}    
\caption{Estimation of  individual patients' flow rate in a breathing cycle for the split ventilator setup, based on noisy measurements of total flow rate and volume at the ventilator.} 
\label{fig:2PatientFlowRate}
\end{center}
\end{figure}

\begin{figure}[t]
\begin{center}
\includegraphics[width=\columnwidth]{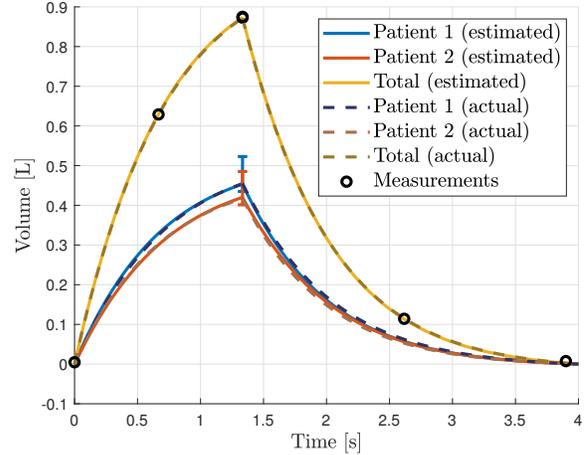}    
\caption{Estimation of  individual patients' flow in a breathing cycle for the split ventilator setup, based only on noisy measurements of total flow rate and volume at the ventilator.} 
\label{fig:2PatientFlow}
\end{center}
\end{figure}

\begin{figure}[t]
\begin{center}
\includegraphics[width=\columnwidth]{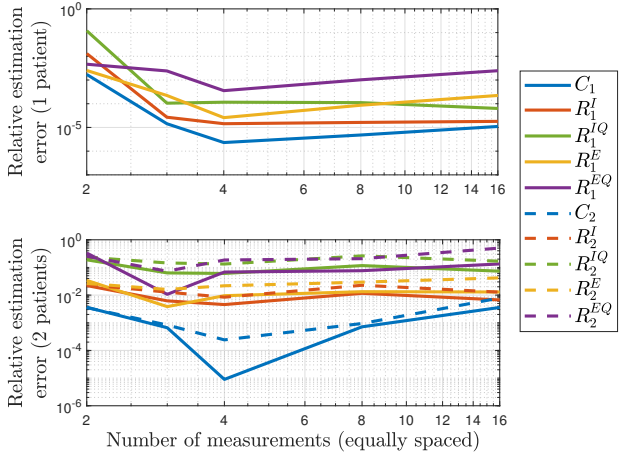}    
\caption{Comparison of achievable estimation errors between the one-patient scenario (with both flow rate and flow measurement) and two-patient scenario (both flow rate and volume measurement of the total values as well as for one splitted branch).} 
\label{fig:EstimationError}
\end{center}
\end{figure}

Numerical results for the estimation case study (Section~\ref{sec:CaseStudy_Estimation}) are shown in Figures~\ref{fig:2PatientFlowRate} and \ref{fig:2PatientFlow} for two patients. The weighted ($S_\nu=I$, $S_w=0.001I$) least-squares cost~\eqref{eqn:ls} is minimized subject to dynamics~\eqref{eqn:dynamics} and measurement constraints \eqref{eqn:measure_flow_rate_ventilator}. We further assume that both patients have reached a limit cycle, and enforce the additional constraints \eqref{eqn:limit_cycle}. In this example, the two-patients are distinct, with true parameter values $C_1 = 0.54$, $C_2=0.49$\,L/cmH$_2$O, and  linear resistances in both paths, $R_1^\varrho= 12.06, \ R_2^\varrho=12.86$\,cmH$_2$O/L/s$, \varrho \in \{I,E\}$.   Simulation data was also generated using the quadratic dynamics detailed in \eqref{eqn:dynamics}, with all resistances $R_p^{IQ}=R_p^{EQ}=2$\,cmH$_2$O/(L/s)$^2, \ p \in\{1,2\}$. We chose $\overline{\nu}=\overline{w}=0.005$. Three noisy measurements with random noise in the range of [$-0.005$, $0.005$] are taken at the ventilator of the combined flow and volume going to both patients in each phase with the measurement times equally spaced. 
Figure~\ref{fig:2PatientFlowRate} shows the actual and estimated air flow rate going to each patient, while Figure~\ref{fig:2PatientFlow} shows the corresponding volumes. In these results we observe that the estimated tidal volume $\Delta Q_E(i_p)$  --- which is a metric of key interest --- is less sensitive to measurement noise than the actual parameter values (see discussion on Figure~\ref{fig:EstimationError}). The error bounds on the tidal volume in Figure~\ref{fig:2PatientFlow} are determined by solving the additional four optimization problems of the form
\begin{equation}
\min \pm C_p\left(v_p(t_0+t_I) – v_p(t_0) \right),
\end{equation}
subject to \eqref{eqn:dynamics}, 
\eqref{eqn:measure_flow_rate_ventilator}.

This gives us the range of admissible tidal volumes that fit the measured data for a prescribed error bound. We can see that with three measurements we can achieve good estimates of the combined and individual flows for both  patients.

\begin{figure}[t]
\begin{center}
\includegraphics[width=\columnwidth]{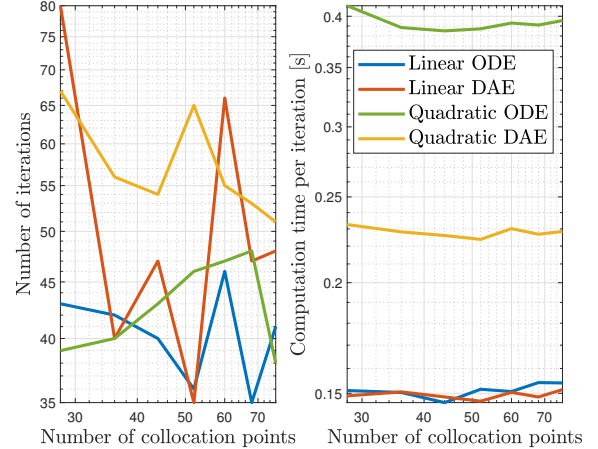}    
\caption{Comparison of computational performance for different model formulations in the two-patient estimation problem, with Hermite-Simpson collocation.} 
\label{fig:2PatientModelFormCompare}
\end{center}
\end{figure}

Next, we demonstrate a number of computational results. To make the solutions reproducible,  we now use perfect measurements with no noise for each phase. However, the problem formulations remain the same with $\overline{\nu}$ and $\overline{w}$ configured to be $0.005$. The two-patient estimation problem~\eqref{eqn:ls} can be solved with different modeling choices. If the quadratic terms in the DAE equations only have minor contributions relative to the linear terms, a linear DAE model can be used to capture the dominating dynamic behaviors. The number of static parameters for the corresponding DOP will therefore be reduced by four, and Figure~\ref{fig:2PatientModelFormCompare} shows that the computational time required per NLP iteration has been significantly decreased. 

For transcription methods that struggle to handle 
DAEs, especially those with high indices, a common practice is to perform index reduction and eventually pose the problem in an ODE only form. However, this process can be tedious for complex systems and often requires the introduction of additional variables, which increase the dimension of the optimization problem considerably. In this example, transforming the quadratic DAE system into an equivalent ODE form requires the number of states to be doubled, and consequently results in considerable computation penalties, as demonstrated in Figure~\ref{fig:2PatientModelFormCompare}. 

In brief, transformation between DAE and ODE for the estimation problem is performed as follows. The inhale phase DAE in \eqref{eqn:dynamics} may be written as
\begin{align*}
    \dot{v}_p(t) &= i_p(t)/C_p, \\ 
     v_p(t) &= -i_p^2(t)R_{pQ}-i(t)R_{pL} + V_I + w(t),
\end{align*}
where $R_{pQ}$ and $R_{pL}$ represent the combined quadratic and linear resistance components respectively. Through application of the chain rule,
\begin{equation*}
    \dot{v}_p(t) = \frac{\mathrm{d}v_p(t)}{\mathrm{d}i_p(t)}\frac{\mathrm{d}i_p(t)}{\mathrm{d}t}, \quad \frac{\mathrm{d}v_p(t)}{\mathrm{d}i_p(t)} = -2i_p(t)R_{pQ} - R_{pL},
\end{equation*}
from which the flow rate dynamics may be written as
\begin{equation*}
    \frac{\mathrm{d}i_p(t)}{\mathrm{d}t} = \frac{-i_p(t)}{C_p(2i_p(t)R_{pQ} + R_{pL})}.
\end{equation*}
We apply the same procedure for the exhale phase. To make the DAE and ODE representations equivalent, we must enforce the original DAE at the beginning of each phase as boundary conditions.

Next, we checked the achievable estimation error under the best measurement paradigm, where both flow rate $i(t)$ and flow volume $Q(t):=\int_{t_0}^t i(\varsigma)d\varsigma$ are measured by the ventilator. In the two-patient scenario additional measurements are available for one of the branches. Figure~\ref{fig:EstimationError} illustrates that, in the case of one patient only, the corresponding patient parameters can be determined very accurately by the proposed estimation scheme. When the parameters of two patients are identified together, although one may not reach the same accuracy level as in one-patient case, the relative estimation errors are still sufficiently small for practical use in most cases.  Moreover, it was found that having more than three measurement points per phase did not lead to any obvious advantages in our tests. 

\subsection{Control}

\begin{figure}[t]
\begin{center}
\includegraphics[width=\columnwidth]{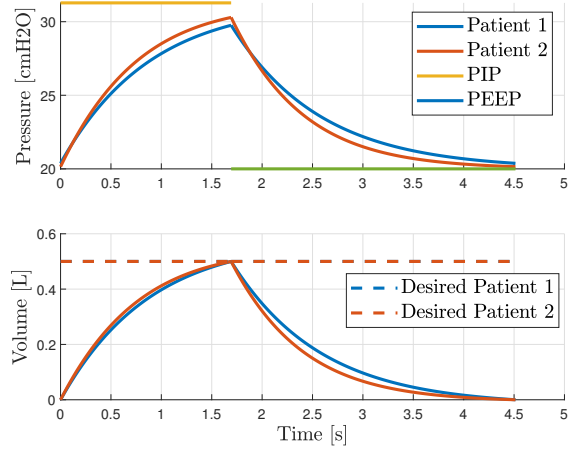}    
\caption{Solution of the ventilator splitting control problem using the minimum energy objective and constant PIP and PEEP pressures in the breathing cycle. This solution requires the adjustable resistance for patient~1 to be set to its maximum for both inhale and exhale cycles, while the resistance for patient~2 should be set to $0$. Total energy usage is 11.3\,J.} 
\label{fig:2PatientControl_ConstVIVE}
\end{center}
\end{figure}

\begin{figure}[t]
\begin{center}
\includegraphics[width=\columnwidth]{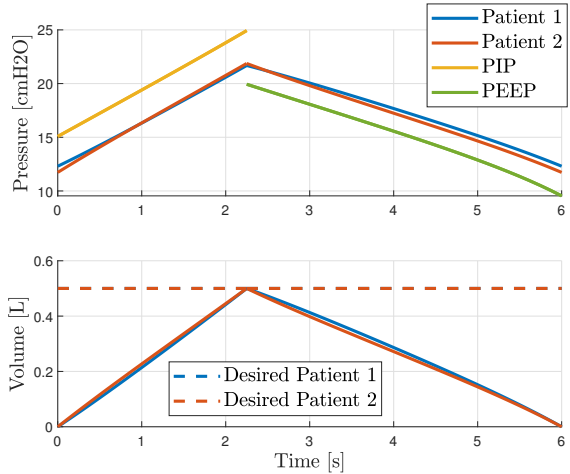}    
\caption{Solution of the ventilator splitting control problem using the minimum energy objective and time-varying PIP and PEEP pressures in the breathing cycle. This solution requires the adjustable resistance to be set to 47\% of its maximum for patient~1 during inhale, and at its maximum during exhale. The adjustable resistance for patient~2 should be configured to $0$ in both phases. Total energy usage is 5.1\,J.} 
\label{fig:2PatientControl_VaryingVIVE}
\end{center}
\end{figure}

Figures~\ref{fig:2PatientControl_ConstVIVE} and \ref{fig:2PatientControl_VaryingVIVE} show the solution of the control problem~\eqref{eqn:Control_Problem} with minimum energy cost for two patients under the respective assumptions of either constant or time-varying PIP and PEEP pressures. 
In both examples the respiratory rate is restricted to $10<f_R<20$ breaths per minute, and the inhale to exhale ratio is in the range $0.4<\rho<0.6$. The bounds on pressures, in cmH$_2$O, are $15<V_I<35$ and $5<V_E<20$.
The patients are parameterized using the $\hat \vartheta(\mathcal{U}_{old},\mathcal{Y})$ values corresponding to the measurements and solution of the estimation problem \eqref{eqn:ls}, shown in Figures~\ref{fig:2PatientFlowRate} and \ref{fig:2PatientFlow}. Importantly, if the resistance  of the adjustable resistances $a_{I,p},a_{E,p}$ are fixed, as is typically the case in one-patient ventilator setups, then we may not be able to achieve arbitrary tidal volumes for all patients. 

However, by including these adjustable resistances in our manipulated variables $\mathcal{U}_{new}$, both patients receive a tidal volume of $0.5$\,L. If the PIP and PEEP pressures are held constant over the inhale and exhale phases, then the solution requires the resistance for patient 1 to be maximized, $a_{\varrho,1}=1$, and the resistance for patient 2 to be minimized, $a_{\varrho,2}=0, \varrho\in\{I,E\}$. As shown in Figure~\ref{fig:2PatientControl_ConstVIVE}, the minimum energy solution results in a PEEP pressure of $V_E=20$\,cmH$_2$O and  PIP pressure of $V_E=31.3$\,cmH$_2$O.

If we now allow for time-varying $V_I$ and $V_E$, then we can achieve the same tidal volume for each patient, with less variation in pressure, as shown in Figure~\ref{fig:2PatientControl_VaryingVIVE}. This may be beneficial for patient comfort and recovery. The corresponding flow appears very similar to a volume-controlled ventilator operation, where the flow rate is typically constant over the inhale phase. This results in an approximately linear increase in the flow delivered to each patient. In this case the adjustable resistances $a_{I,1}=0.47, a_{E,1}=1$, while $a_{I,2} = a_{E,2}=0$. Allowing for time-varying pressures results in less than $50\%$ of the energy used when $V_I$ and~$V_E$ were held constant over each phase.

\begin{figure}[t]
\begin{center}
\includegraphics[width=\columnwidth]{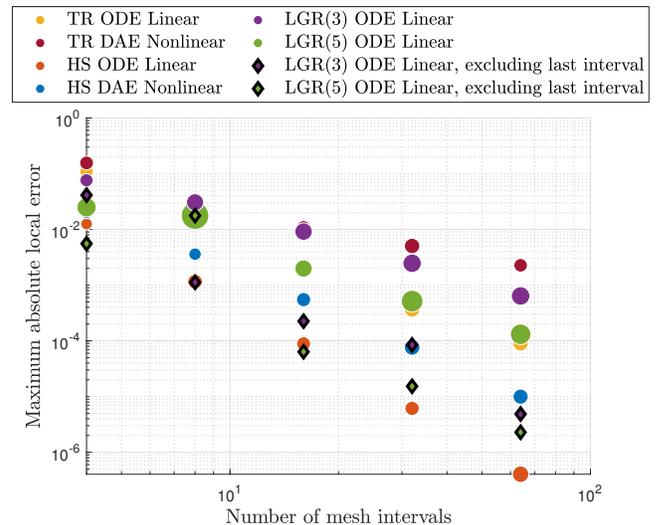}    
\caption{Comparison of computational performance and solution accuracy for different model formulations, discretization methods and mesh sizes, for the ventilator splitting control problem solved using direct collocation and without warm-starting. The radius of the dot is proportional to the solution computation time, ranging from 5.8\,s (corresponding to HS nonlinear DAE  with 4 intervals) to 31.7\,s (corresponding to a failed convergence case of LGR(5) with 8 intervals). TR: Trapezoidal; HS: Hermite-Simpson; LGR($p$): LGR with degree $p$ polynomial.} 
\label{fig:2PatientControlMaxErrorComputation}
\end{center}
\end{figure}

Figure~\ref{fig:2PatientControlMaxErrorComputation} demonstrates the relationship between the mesh density and the achievable error levels using direct collocation transcription with various discretization schemes. 

The error measure used is the maximum absolute local error, introduced in Section~\ref{subsec:Error_And_MeshRefinement}, considering both the ODE and DAE (equality constraint) residuals and inequality constraint violations. It can be observed that higher-order discretization methods generally lead to smaller errors, compared to lower-order methods on the same mesh, at the cost of being computationally more expensive. One exception is when comparing the  trapezoidal and Hermite-Simpson methods to the LGR method. For LGR, the last node of the mesh is not a collocation point. As a result, additional challenges arises when enforcing the periodic boundary conditions, with the absolute local error for the last interval a few magnitudes larger than elsewhere. Therefore, the maximum absolute local errors for LGR look particularly worse when compared to the other methods, especially considering the high degree of polynomials employed. This observation highlights that collocation methods with non-collocating end-point constraints may be unsuitable for particular problem formulations. Increasing the number of mesh intervals, i.e.\ making the mesh denser, will generally result in improved accuracy, but at the cost of increased computational effort.

It is not difficult, however, to find cases that do not follow this general trend. Sometimes the computation can be faster for the dense mesh and vice versa. A scheme which initializes the solve on a dense grid using a lower accuracy solution obtained from a coarse grid can often lead to reduced computation compared to solving on the same dense grid with a poor initialization. As a result, mesh refinement strategies have become a crucial aspect in designing efficient numerical methods to solve large-scale DOPs.

\begin{figure}[t]
\begin{center}
\includegraphics[width=\columnwidth]{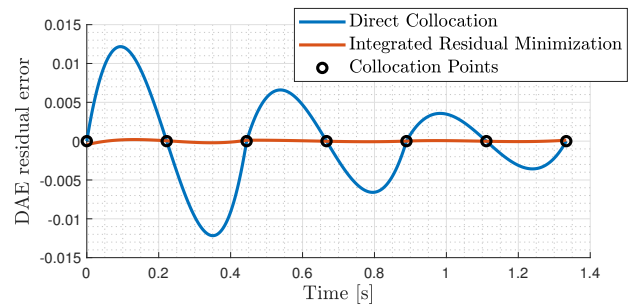}    
\caption{Comparison of DAE residual errors along the trajectory of the inhale phase of the ventilator splitting control problem, for two different transcription methods. Both solutions are piecewise cubic with 3 mesh intervals.} 
\label{fig:DAEResidualError}
\end{center}
\end{figure}

Figure~\ref{fig:DAEResidualError} shows the DAE residual errors corresponding to the solution trajectory for a single inhale phase computed using both direct collocation and direct integrated residual minimization on an extremely coarse grid. Figure~\ref{fig:DAEResidualError} illustrates the fundamental differences between these transcription methods.
It is clear that for direct collocation, despite forcing the DAE residuals to be zero at collocation points, large errors can still occur in-between them. A method that minimizes the total residual error integrated along the whole trajectory can yield solutions of much higher accuracy for the same coarse mesh. This can be beneficial for many embedded applications where computational resources are limited, preventing the use of a dense mesh. 

\section{Conclusions}

\subsection{Direct Transcription Methods}

In the first part we presented an overview of different direct methods for transforming and discretizing a nonlinear, continuous-time dynamic optimization problem into a structured nonlinear program, which can be efficiently solved using state-of-the-art numerical methods. We briefly introduced error analysis and mesh refinement schemes, which are important in achieving a desired accuracy to the dynamic optimization problem, and discussed why appropriately-defined integrals of the differential equation residual are suitable for error analysis.

Collocation is the easiest transcription method to implement and  allows for a large class of parameterizations for the solution trajectories, hence why it has been very successful in a number of applications.  However, it can be difficult to get collocation methods to work with DAEs and certain classes of problems, due to requiring that the residual be exactly zero  at collocation points. As a consequence, there is a complex interplay between the choice of collocation points and parameterizations used. 

Integrated residual methods, though not yet as popular as collocation methods, allows for a more straightforward handling of challenging problems and can result in smaller residual errors, compared to collocation methods. However, the implementation of integrated residual methods is slightly more involved than for collocation methods. Further research is therefore needed to make it easier for the wider deployment of integrated residual methods.

Runge-Kutta methods might be the most well-known class of methods for solving differential equations. It might therefore be tempting to a dynamic optimization novice to start with this class of methods. However, if the chosen Runge-Kutta method is not equivalent to a polynomial collocation method, then it can be difficult to define an appropriate  parameterization and  error analysis scheme, while ensuring that the implementation is as efficient as a polynomial collocation method of similar order. This is why, in dynamic optimization problems, polynomial collocation methods are often preferred over other Runge Kutta methods. Having said this, it is always worth trying out different Runge-Kutta methods if  collocation methods fail.

Shooting and simultaneous transcription methods have their respective pros and cons. Multiple shooting only parameterizes free variables inside mesh intervals, resulting in smaller NLP problems. Also, it is easier to make code for shooting methods modular. The ability to easily interface to state-of-the art differential equation solvers  explains the success of shooting methods. On the other hand, simultaneous transcription methods parameterize the state variables in addition to the free variables, which  generally results in larger, but sparser NLP problems. Depending on the problem, this larger NLP may be solved more efficiently with sparsity-exploiting solvers, compared to the smaller NLP of a shooting method. Simultaneous methods also avoid the dependency on separate ODE/DAE solvers, as in shooting  methods. Hence, simultaneous methods could result in simpler, standalone schemes. 

The list of topics discussed here is nowhere complete. For example, we omitted a detailed discussion on robust and stochastic dynamic optimization problems, which has seen a growth in interest over the last few years. We also did not discuss multi-phase problems, where state trajectories are allowed to be discontinuous, but hopefully the interested reader now has the framework in which to explore this. We also chose not to discuss tailored numerical methods for solving the resulting NLPs, but we hope that some of the references we provided in this paper will serve as a good starting point on this topic.

\subsection{Dual-patient Ventilation}

We explored the application of some of the discussed methods to a simulation case study on the system identification, estimation and control of ventilators being split between multiple patients.
We showed that, by solving a suitably-defined dynamic optimization problem, it is possible to accurately reconstruct the flow being delivered to each patient, even if we do not have individual measurements of the flow to each patient. A direct byproduct of solving this problem is an estimate of the resistance and compliance components used to model the flow for each patient. With these estimates, we can modify the resistance in the inhale and exhale path of each patient, in order to ensure that the required tidal volume is   delivered to each patient.  

We would like to emphasise that ventilating multiple patients from a single ventilator is currently untested. The results are of a simulation study and do not include any experimental component. The work presented here does not change the current clinical guidelines. The main aim of this paper was to serve as a tutorial to dynamic optimization.  We hope, however, that the application of these advanced numerical methods to split ventilation has shed some scientific light on a very challenging, topical problem.

\balance

\appendices
 
\section{Lagrange Polynomials}
\label{app:Lagrange}

Suppose we are given $n_v+1$ values $v_j\in\mathbb{R}$ at distinct nodes $t_j$, $j\in\mathbb{J}:=\{0,\ldots,n_v\}$. The lowest degree polynomial $t\mapsto p(t)$ that interpolates the values  at these points, i.e.
\[
p(t_j) = v_j,\ \forall j \in\mathbb{J},
\]
is given by the following polynomial in \emph{Lagrange form}:
\[
    p(t):=\sum_{j\in\mathbb{J}}  v_j \mathcal{L}_j(t) ,
\]
where the \emph{Lagrange polynomial} $\mathcal{L}_j$ corresponding to the node~$t_j$ is given by
\[
\mathcal{L}_j(t) :=  \prod_{k\in\mathbb{J},j\neq k} \frac{t-t_k}{t_j-t_k},\ \forall j \in\mathbb{J}.
\]
Evaluating $p(t)$ or $\dot{p}(t)$, given  the set $\{(t_j,p_j)\mid j\in\mathbb{J}\}$, can be done efficiently and in a numerically stable manner via the \emph{barycentric formula}~\cite{doi:10.1137/S0036144502417715}.

\bibliographystyle{IEEEtran}
\bibliography{CDC2020tutorial,medical}

\begin{thebibliography}{10}
\providecommand{\url}[1]{#1}
\csname url@rmstyle\endcsname
\providecommand{\newblock}{\relax}
\providecommand{\bibinfo}[2]{#2}
\providecommand\BIBentrySTDinterwordspacing{\spaceskip=0pt\relax}
\providecommand\BIBentryALTinterwordstretchfactor{4}
\providecommand\BIBentryALTinterwordspacing{\spaceskip=\fontdimen2\font plus
\BIBentryALTinterwordstretchfactor\fontdimen3\font minus
  \fontdimen4\font\relax}
\providecommand\BIBforeignlanguage[2]{{%
\expandafter\ifx\csname l@#1\endcsname\relax
\typeout{** WARNING: IEEEtran.bst: No hyphenation pattern has been}%
\typeout{** loaded for the language `#1'. Using the pattern for}%
\typeout{** the default language instead.}%
\else
\language=\csname l@#1\endcsname
\fi
#2}}

\bibitem{betts2010practical}
J.~Betts, \emph{Practical Methods for Optimal Control and Estimation Using
  Nonlinear Programming: Second Edition}, ser. Advances in Design and
  Control.\hskip 1em plus 0.5em minus 0.4em\relax Society for Industrial and
  Applied Mathematics, 2010.

\bibitem{RMD2nd}
J.~B. Rawlings, D.~Q. Mayne, and M.~Diehl, \emph{Model Predictive Control:
  Theory, Computation, and Design}, 2nd~ed.\hskip 1em plus 0.5em minus
  0.4em\relax Nob Hill Publishing, Santa Barbara CA, USA, 2019.

\bibitem{HandbookMPC}
S.~V. Rakovic and W.~S. Levine, Eds., \emph{Handbook of Model Predictive
  Control}.\hskip 1em plus 0.5em minus 0.4em\relax Birk\"auser, 2019.

\bibitem{doi:10.1098/rsos.200585}
\BIBentryALTinterwordspacing
J.~A. Solís-Lemus, E.~Costar, D.~Doorly, E.~C. Kerrigan, C.~H. Kennedy,
  F.~Tait, S.~Niederer, P.~E. Vincent, and S.~E. Williams, ``A simulated single
  ventilator/dual patient ventilation strategy for acute respiratory distress
  syndrome during the {COVID}-19 pandemic,'' \emph{Royal Society Open Science},
  vol.~7, no.~8, p. 200585, 2020. [Online]. Available:
  \url{https://royalsocietypublishing.org/doi/abs/10.1098/rsos.200585}
\BIBentrySTDinterwordspacing

\bibitem{FERREAU201713194}
\BIBentryALTinterwordspacing
H.~Ferreau, S.~Alm{\'e}r, R.~Verschueren, M.~Diehl, D.~Frick, A.~Domahidi,
  J.~Jerez, G.~Stathopoulos, and C.~Jones, ``Embedded optimization methods for
  industrial automatic control,'' \emph{IFAC-PapersOnLine}, vol.~50, no.~1, pp.
  13\,194 -- 13\,209, 2017, 20th IFAC World Congress. [Online]. Available:
  \url{http://www.sciencedirect.com/science/article/pii/S2405896317325764}
\BIBentrySTDinterwordspacing

\bibitem{iet:/content/journals/10.1049/iet-cta.2019.0168}
\BIBentryALTinterwordspacing
A.~Beghi, ``\BIBforeignlanguage{English}{Efficient move blocking strategy for
  multiple shooting-based non-linear model predictive control},''
  \emph{\BIBforeignlanguage{English}{IET Control Theory \& Applications}},
  vol.~14, pp. 343--351(8), January 2020. [Online]. Available:
  \url{https://digital-library.theiet.org/content/journals/10.1049/iet-cta.2019.0168}
\BIBentrySTDinterwordspacing

\bibitem{doi:10.1137/S0036144502417715}
\BIBentryALTinterwordspacing
J.-P. Berrut and L.~N. Trefethen, ``Barycentric {L}agrange interpolation,''
  \emph{SIAM Review}, vol.~46, no.~3, pp. 501--517, 2004. [Online]. Available:
  \url{https://doi.org/10.1137/S0036144502417715}
\BIBentrySTDinterwordspacing

\bibitem{liu2017adaptive}
F.~Liu, W.~W. Hager, and A.~V. Rao, ``Adaptive mesh refinement method for
  optimal control using decay rates of legendre polynomial coefficients,''
  \emph{IEEE Transactions on Control Systems Technology}, vol.~26, no.~4, pp.
  1475--1483, 2017.

\bibitem{ECH_CSL:2020}
Y.~{Nie} and E.~C. {Kerrigan}, ``External constraint handling for solving
  optimal control problems with simultaneous approaches and interior point
  methods,'' \emph{IEEE Control Systems Letters}, vol.~4, no.~1, pp. 7--12,
  2020.

\bibitem{polak:2009}
H.~{Chung}, E.~{Polak}, and S.~{Sastry}, ``An external active-set strategy for
  solving optimal control problems,'' \emph{IEEE Transactions on Automatic
  Control}, vol.~54, no.~5, pp. 1129--1133, 2009.

\bibitem{FP2019}
F.~A. C.~C. {Fontes} and L.~T. {Paiva}, ``Guaranteed constraint satisfaction in
  continuous-time control problems,'' \emph{IEEE Control Systems Letters},
  vol.~3, no.~1, pp. 13--18, 2019.

\bibitem{DR2ns}
P.~J. Davis and P.~Rabinowitz, \emph{Methods of Numerical Integration},
  2nd~ed.\hskip 1em plus 0.5em minus 0.4em\relax Dover Publications, USA, 1984.

\bibitem{kelly2017}
\BIBentryALTinterwordspacing
M.~Kelly, ``An introduction to trajectory optimization: How to do your own
  direct collocation,'' \emph{SIAM Rev.}, vol.~59, no.~4, pp. 849--904, 2017.
  [Online]. Available: \url{https://doi.org/10.1137/16M1062569}
\BIBentrySTDinterwordspacing

\bibitem{garg2009overview}
D.~Garg, M.~Patterson, W.~Hager, A.~Rao, D.~Benson, and G.~Huntington, ``An
  overview of three pseudospectral methods for the numerical solution of
  optimal control problems,'' \emph{Advances in the Astronautical Sciences},
  vol. 135, no.~1, pp. 475--487, 2009.

\bibitem{fahroo2008advances}
F.~Fahroo and I.~Ross, ``Advances in pseudospectral methods for optimal
  control,'' in \emph{AIAA guidance, navigation and control conference and
  exhibit}, 2008, p. 7309.

\bibitem{NieKerriganCSL:2020}
Y.~{Nie} and E.~C. {Kerrigan}, ``Efficient and more accurate representation of
  solution trajectories in numerical optimal control,'' \emph{IEEE Control
  Systems Letters}, vol.~4, no.~1, pp. 61--66, Jan 2020.

\bibitem{NeuenhofenKerrigan:2018}
M.~P. Neuenhofen and E.~C. Kerrigan, ``Dynamic optimization with convergence
  guarantees,'' \emph{arXiv preprint arXiv:1810.04059}, 2018.

\bibitem{NeuenhofenKerrigan:CDC2020}
------, ``An integral penalty-barrier direct transcription method for optimal
  control,'' in \emph{Proc.\ 59th IEEE Conference on Decision and Control},
  2020.

\bibitem{kunkelMehrmann2006}
P.~Kunkel and V.~Mehrmann, \emph{Differential-Algebraic Equations: Analysis and
  Numerical Solution}.\hskip 1em plus 0.5em minus 0.4em\relax European
  Mathematical Society, 2006.

\bibitem{IPOPT}
A.~W{\"a}chter and L.~T. Biegler, ``On the implementation of an interior-point
  filter line-search algorithm for large-scale nonlinear programming,''
  \emph{Mathematical programming}, vol. 106, no.~1, pp. 25--57, 2006.

\bibitem{SNOPT}
P.~E. Gill, W.~Murray, and M.~A. Saunders, ``{SNOPT}: An {SQP} algorithm for
  large-scale constrained optimization,'' \emph{SIAM review}, vol.~47, no.~1,
  pp. 99--131, 2005.

\bibitem{WORHP}
C.~Büskens and D.~Wassel, ``The {ESA} {NLP} solver {WORHP},'' in
  \emph{Modeling and Optimization in Space Engineering}, G.~Fasano and J.~D.
  Pintér, Eds.\hskip 1em plus 0.5em minus 0.4em\relax Springer New York, 2013,
  vol.~73, pp. 85--110.

\bibitem{NOMAD}
S.~Le~Digabel, ``Algorithm 909: Nomad: Nonlinear optimization with the mads
  algorithm,'' \emph{ACM Transactions on Mathematical Software (TOMS)},
  vol.~37, no.~4, pp. 1--15, 2011.

\bibitem{NWKSB}
\BIBentryALTinterwordspacing
B.~L. Nicholson, W.~Wan, S.~Kameswaran, and L.~T. Biegler, ``Parallel cyclic
  reduction strategies for linear systems that arise in dynamic optimization
  problems,'' \emph{Computational Optimization and Applications}, vol.~70,
  no.~2, pp. 321--350, 2018. [Online]. Available:
  \url{https://doi.org/10.1007/s10589-018-0001-7}
\BIBentrySTDinterwordspacing

\bibitem{doi:10.1002/jmv.25722}
\BIBentryALTinterwordspacing
P.~Sun, X.~Lu, C.~Xu, W.~Sun, and B.~Pan, ``Understanding of {COVID-19} based
  on current evidence,'' \emph{Journal of Medical Virology}, vol.~92, no.~6,
  pp. 548--551, 2020. [Online]. Available:
  \url{https://onlinelibrary.wiley.com/doi/abs/10.1002/jmv.25722}
\BIBentrySTDinterwordspacing

\bibitem{WHO}
\BIBentryALTinterwordspacing
{World Health Organization}. (2020, May) Coronavirus disease ({COVID-19})
  pandemic. [Online]. Available:
  \url{https://www.who.int/emergencies/diseases/novel-coronavirus-2019}
\BIBentrySTDinterwordspacing

\bibitem{CCDCP}
\BIBentryALTinterwordspacing
{Chinese Center for Disease Control and Prevention}, ``Distribution of new
  coronavirus pneumonia.'' [Online]. Available:
  \url{http://2019ncov.chinacdc.cn/2019-nCoV/}
\BIBentrySTDinterwordspacing

\bibitem{JHCHS}
\BIBentryALTinterwordspacing
{Johns Hopkins Center for Health Security}. (2020, April) Ventilator
  stockpiling and availability in the {US}. [Online]. Available:
  \url{https://www.centerforhealthsecurity.org/resources/COVID-19/COVID-19-fact-sheets/200214-VentilatorAvailability-factsheet.pdf}
\BIBentrySTDinterwordspacing

\bibitem{doi:10.1056/NEJMp2006141}
\BIBentryALTinterwordspacing
M.~L. Ranney, V.~Griffeth, and A.~K. Jha, ``Critical supply shortages --- the
  need for ventilators and personal protective equipment during the {Covid-19}
  pandemic,'' \emph{New England Journal of Medicine}, vol. 382, no.~18, p. e41,
  2020. [Online]. Available: \url{https://doi.org/10.1056/NEJMp2006141}
\BIBentrySTDinterwordspacing

\bibitem{ASA}
\BIBentryALTinterwordspacing
{American Society of Anesthesiologists}, March 2020. [Online]. Available:
  \url{https://www.asahq.org/about-asa/newsroom/news-releases/2020/03/joint-statement-on-multiple-patients-per-ventilator}
\BIBentrySTDinterwordspacing

\bibitem{tronstad2020splitting}
C.~Tronstad, T.~Martinsen, M.~Olsen, L.~Rosseland, F.~Pettersen,
  {\O}.~Martinsen, J.~H{\o}getveit, and H.~Kalv{\o}y, ``Splitting one
  ventilator for multiple patients--a technical assessment,'' \emph{arXiv
  preprint arXiv:2003.12349}, 2020.

\bibitem{campbell1963electrical}
D.~Campbell and J.~Brown, ``The electrical analogue of lung,'' \emph{British
  Journal of Anaesthesia}, vol.~35, no.~11, pp. 684--692, 1963.

\bibitem{Plummer2020.04.12.20062497}
\BIBentryALTinterwordspacing
A.~R. Plummer, J.~L. {du Bois}, J.~M. Flynn, J.~Roesner, S.~M. Lee, P.~Magee,
  M.~Thornton, A.~Padkin, and H.~S. Gill, ``A simple method to estimate flow
  restriction for dual ventilation of dissimilar patients: The {BathRC}
  model,'' \emph{medRxiv}, 2020. [Online]. Available:
  \url{https://doi.org/10.1101/2020.04.12.20062497}
\BIBentrySTDinterwordspacing

\bibitem{doi:10.1164/rccm.201612-2495CI}
\BIBentryALTinterwordspacing
W.~R. Henderson, L.~Chen, M.~B.~P. Amato, and L.~J. Brochard, ``Respiratory
  mechanics in acute respiratory distress syndrome,'' \emph{American Journal of
  Respiratory and Critical Care Medicine}, vol. 196, no.~7, pp. 822--833, 2017,
  pMID: 28306327. [Online]. Available:
  \url{https://doi.org/10.1164/rccm.201612-2495CI}
\BIBentrySTDinterwordspacing

\bibitem{nie2018iclocs2}
Y.~Nie, O.~Faqir, and E.~C. Kerrigan, ``{ICLOCS2}: Try this optimal control
  problem solver before you try the rest,'' in \emph{2018 UKACC 12th
  International Conference on Control (CONTROL)}.\hskip 1em plus 0.5em minus
  0.4em\relax IEEE, 2018, pp. 336--336.

\end{thebibliography}

\end{document}